\documentclass[11pt]{article}
\usepackage{amsmath}
\usepackage{amssymb}
\usepackage{amsthm}
\usepackage{latexsym}
\usepackage{graphicx}
\usepackage{hyperref}

\newtheorem{thm}{Theorem}[section]

\newtheorem{lem}[thm]{Lemma}
\newtheorem{prop}[thm]{Proposition}
\newtheorem{defn}[thm]{Definition}

\newcommand{\mb}{\mathbf}
\newcommand{\mh}{\mathbb}
\newcommand{\mr}{\mathrm}
\newcommand{\mc}{\mathcal}
\newcommand{\mf}{\mathfrak}

\newcommand{\ts}{\textstyle}

\newcommand{\Z}{\mathbb Z}
\newcommand{\Q}{\mathbb Q}
\newcommand{\R}{\mathbb R}
\newcommand{\C}{\mathbb C}

\newcommand{\es}{\emptyset}

\newcommand{\pch}{periodic cyclic homology}

\newcommand{\inp}[2]{\langle #1 \,,\, #2 \rangle}

\newcommand{\hexagon}[6]{\begin{array}{ccccc}
  #1 & \to & #2 & \to & #3 \\
  \uparrow & & & & \downarrow \\
  #4 & \leftarrow & #5 & \leftarrow & #6
  \end{array}}

\begin{document}

\begin{center}
\LARGE \textbf{Homology of graded Hecke algebras}\\[4mm]
Maarten Solleveld \\[15mm]
\normalsize Mathematisches Institut,
Georg-August-Universit\"at G\"ottingen\\
Bunsenstra\ss e 3-5, 37073 G\"ottingen, Germany\\
email: Maarten.Solleveld@mathematik.uni-goettingen.de \\
October 2009\\[14mm]
\end{center}
\textbf{Abstract.}\\
Let $\mh H$ be a graded Hecke algebra with complex deformation parameters and 
Weyl group $W$. We show that the Hochschild, cyclic and periodic cyclic homologies 
of $\mh H$ are all independent of the parameters, and compute them explicitly.
We use this to prove that, if the deformation parameters are real, the collection of 
irreducible tempered $\mh H$-modules with real central character forms a 
$\Q$-basis of the representation ring of $W$.

Our method involves a new interpretation of the \pch \ of finite type algebras,
in terms of the cohomology of a sheaf over the underlying complex affine variety.

\tableofcontents

\newpage

\section*{Introduction}
\addcontentsline{toc}{section}{Introduction}

Let $\mf t^*$ be a complex vector space containing a root system $R$ with Weyl group $W$.
Then $W$ acts on the symmetric algebra $S (\mf t^* )$ of $\mf t^*$, so one can construct
the crossed product algebra $W \ltimes S(\mf t^*)$. Graded Hecke algebras are deformations of 
$W \ltimes S (\mf t^*)$, depending on several parameters $k_\alpha \in \mh C$. 

Lusztig \cite{Lus1,Lus2} showed that graded Hecke algebras play an important role in the
representation theory of affine Hecke algebras and of simple $p$-adic groups. For this
reason it pays to understand the representations of graded Hecke algebras, and in
particular to determine their irreducible representations.

We will approach this goal via the noncommutative geometry of these algebras. Every
graded Hecke algebra $\mh H$ is endowed with a natural filtration, whose associated graded
algebra is $W \ltimes S (\mf t^* )$. This gives rise to spectral sequences converging
to various homologies of $\mh H$, and allows us to prove that:

\begin{thm}\label{thm:0.1} (See Theorems \ref{thm:4.3} and \ref{thm:4.4}.)\\
The Hochschild and (periodic) cyclic homologies of $\mh H$ and $W \ltimes S (\mf t^* )$
are isomorphic:
\[
\begin{array}{ccc}
HH_n (\mh H ) & \cong & HH_n (W \ltimes S(\mf t^* )) , \\
HC_n (\mh H ) & \cong & HC_n (W \ltimes S(\mf t^* )) , \\
HP_n (\mh H ) & \cong & HP_n (W \ltimes S(\mf t^* )) . \\
\end{array}
\]
Moreover the inclusion $\C [W] \to \mh H$ induces an isomorphism on $HP_*$.
\end{thm}

However, the representation theoretic content of this result is not so obvious, for instance 
because the bijection on $HH_*$ is not an isomorphism of $S( \mf t^* )^W$-modules.
To clarify this we introduce a new interpretation of the \pch \ of an algebra. It is
not specific for graded Hecke algebras, in fact it applies to all finite type algebras.
Let $X$ be a complex affine algebraic variety. It is known \cite{KNS} that 
\begin{equation}\label{eq:0.1}
HP_n (\mc O (X)) \cong {\ts \bigoplus_{m \in \Z}} H^{n + 2m} (X^{an} ; \C ) ,
\end{equation}
where $X^{an}$ is $X$ with the analytic topology. Now let $A$ be a finite type 
$\mc O (X)$-algebra. To it we associate a sheaf $S_A$ over $X^{an}$, which encodes the
space of irreducible $A$-modules. We will establish an isomorphism 
\begin{equation}\label{eq:0.2}
HP_n (A) \cong {\ts \bigoplus_{m \in \Z}} \check H^{n + 2m} (X^{an} ; S_A ) ,
\end{equation}
where $\check H$ stands for the \v Cech cohomology of a sheaf.
The main advantage of \eqref{eq:0.2} is that the extensive machinery of sheaf cohomology
becomes available for the calculation of \pch .

For graded Hecke algebras we can really put \eqref{eq:0.2} into effect, provided that the root
system $R$ is crystallographic and the parameters $k_\alpha$ are real. Probably these
assumptions are not necessary for the final results, but we do use them in our proofs.

The representations of such graded Hecke algebras naturally come in series. Every such 
series is induced from a discrete series representation of a parabolic subalgebra and
parametrized by a subspace of $\mf t$. The individual representations in such a series 
can be equivalent, and they may be reducible. All this is described by the theory of 
intertwining operators, which the author discussed in the prequel \cite{Sol2} to the
present paper. The upshot is that these series partition the primitive ideal spectrum
of $\mh H$ in disjoint subsets, and that every single such subset is in first approximation
a vector space modulo a finite group.

This relatively simple shape allows us to compute the cohomology of the sheaf $S_{\mh H}$ 
from \eqref{eq:0.2}. As could be expected, there is no cohomology in strictly positive
degrees. In degree zero we find that $\check H^0 (\mf t / W ; S_{\mh H})$ is naturally
isomorphic to the complex vector space whose basis is formed by the irreducible tempered
$\mh H$-modules with real central character (these notions are defined in Section 5).

Thus we have two different ways to determine $HP_* (\mh H)$:
\begin{itemize}
\item via a spectral sequence coming from the filtration of $\mh H$ by degree,
\item via the cohomology of a sheaf encoding the primitive ideal spectrum of $\mh H$.
\end{itemize}
The first method relates $HP_0 (\mh H)$ to $HP_0 \big( \C [W] \big)$, and hence to the 
irreducible $W$-representations. The second strategy results in a comparison between 
$HP_0 (\mh H)$ and the collection $\mr{Irr}_0 (\mh H)$ of irreducible tempered $\mh H$-modules 
with real central character (modulo equivalence). 
Combining both strategies, we arrive at the main result of the paper:

\begin{thm}\label{thm:0.2} (See Theorem \ref{thm:6.5}.c.) \\
Let $R$ be a crystallographic root system with Weyl group $W$ and let 
$k_\alpha \;(\alpha \in R)$ be real parameters. The representations $\{ V |_W : 
V \in \mr{Irr}_0 (\mh H) \}$ form a $\Q$-basis of the representation ring $R(W) \otimes_\Z \Q$.
\end{thm}

Hitherto this result was only known for ``geometric'' graded Hecke algebras, for which Ciubotaru 
\cite[Corollary 3.6]{Ciu} derived it from deep results of Lusztig. Here ``geometric'' means that 
these are the graded Hecke algebras that arise as the equivariant homology of specific algebraic 
varieties, as in \cite{Lus1}.

Moreover in these cases the expression of $\{ V |_W : V \in \mr{Irr}_0 (\mh H) \}$ with respect to
the natural basis of $R(W)$ (suitably ordered) is a unitriangular matrix with integral coefficients. 
As Ciubotaru remarks, unitriangularity gives rise to a generalized Springer correspondence,
via \cite{Lus0}. Therefore it would be quite interesting to know whether this unitriangularity
remains valid in our more general setting. Unfortunately the methods in this paper are unsuitable
to detect it, since the complex coefficients of \pch \: destroy the torsion information. 

Theorem \ref{thm:0.2} can be used to study several problems:
\begin{itemize}
\item determination of the unitary dual of graded Hecke algebras. Indeed, Theorem \ref{thm:0.2}
(for geometric Hecke algebras) plays an important role in \cite{Ciu}.
\end{itemize}
Lusztig's reduction theorems \cite{Lus2} provide a concrete connection between representations of 
$\mh H$ and of related affine Hecke algebras. This opens the road to
\begin{itemize}
\item the analogue of Theorem \ref{thm:0.1} for affine Hecke algebras.
\end{itemize}
Such an analogue would be useful to classify all irreducible representations of affine Hecke algebras, 
especially those with unequal parameters. Pursuing this line of thought even further, to the smooth 
representation theory of $p$-adic groups, could provide
\begin{itemize}
\item a promising tool to attack many cases of a conjecture of Aubert, Baum and Plymen \cite{ABP,ABP2},
which describes Bernstein components in the smooth dual of a reductive $p$-adic group.
\end{itemize}
The author intends to study the last two applications in a forthcoming paper.

We conclude this introduction with a description of the various sections.
The first Section recalls some general results on crossed products of commutative algebras with
finite groups, in particular their representation theory and their homology. We state the definitions
of graded Hecke algebras in Section 2. As Lusztig's reduction theorems involve groups of diagram
automorphisms of $R$, we incorporate such groups in our algebras right from the start. This changes 
the notations somewhat, but does not make the proofs really more complicated. Section 3 is dedicated 
to proving Theorem \ref{thm:0.1}. In the fourth Section we show how to relate the \pch \ of a finite
type algebra to the cohomology of a sheaf, as in \eqref{eq:0.2}. Moreover we provide some examples
of such sheaves. In Section 5 we recall the relevant results from \cite{Sol2}, which deal mostly with 
parabolically induced representations of graded Hecke algebras. In the final Section 6 we
compute the cohomology of the sheaf $S_{\mh H}$, which leads to Theorem \ref{thm:0.2}.
\vspace{4mm}

\section{Crossed product algebras}
\label{sec:1}

We start with a discussion of a rather general kind of algebras,
obtained from an action of a finite group on a commutative algebra.
A lot is known about such crossed product algebras, and here we 
recall the results that we will use later on.

Let $G$ be any group and $A$ be an algebra over a field $\mh F$. Let 
$\beta : G \to \mr{Aut}(A)$ be an action of $G$ on $A$ by algebra automorphisms. 
We endow the vector space $A \otimes_{\mh F} \mh F G$ with the multiplication
\begin{equation}\label{eq:1.5}
(a \otimes g) \cdot \big( a' \otimes g' \big) = a \, \beta_g (a') \otimes g g'  
\qquad a,a' \in A , g,g' \in G.
\end{equation}
This defines an associative $\mh F$-algebra, denoted $A \rtimes G$ or $G \ltimes A$ 
and called the crossed product of $A$ and $G$.

We specialize this in the direction of the algebras that are of most interest to us.
Let $X$ be a topological space and suppose that $A$ is a subalgebra of $C (X;\C)$
whose maximal ideal spectrum is precisely $X$. Furthermore we suppose that $G$ 
is finite and acts on $X$ by homeomorphisms, such that the induced action on 
$C (X;\C)$ preserves $A$. Let $\C_x$ be the onedimensional $A$-module with character 
$x \in X$ and write 
\begin{align*}
& G_x := \{ g \in G : g(x) = x \} , \\ 
& I_x := \mr{Ind}_A^{A \rtimes G} \C_x \,.
\end{align*}
The representation theory of such algebras is not complicated, and can be 
obtained from classical results that go back to Frobenius and Clifford, see \cite{Cli}
or \cite[Theorem 11.1]{CuRe}.

\begin{thm}\label{thm:1.1}
\begin{description}
\item[a)] $I_x \cong I_{x'}$ if and only if $G x = G x'$.
\item[b)] $I_x \cong \mr{Ind}_{A \rtimes G_x}^{A \rtimes G} \big( \C [G_x] \big)$,
where $A$ acts on $\C [G_x]$ through the evaluation at $x$.
\item[c)] $I_x$ is completely reducible.
\item[d)] Every irreducible $A \rtimes G$-module is a direct summand of some $I_x$.
\item[e)] The number of inequivalent irreducible constituents of $I_x$ equals the
number of conjugacy classes in $G_x$.
\end{description}
\end{thm}
\emph{Proof.}
a) As an $A$-module
\[
{\ts I_x = \bigoplus_{g \in G} \C_{gx} = \bigoplus_{g \in G} \C v_{gx} \,,}
\]
and $G$ acts by $g' \cdot v_{gx} = v_{g' g x}$.
Similarly $I_{x'} = \bigoplus_{g \in G} \C \tilde v_{g x'}$.  If $Gx = Gx'$ then
\[
I_x \to I_{x'} : v_{gx} \mapsto \tilde v_{gx} 
\]
is an isomorphism of $A \rtimes G$-modules. On the other hand, if $G x \neq G x'$, then $I_{x'}$
contains no $A$-weight spaces with weight $x$, so $I_{x'} \not\cong I_x$. \\
b)
\[
\mr{Ind}_{A \rtimes G_x}^{A \rtimes G} \big( \C [G_x] \big) = 
(A \rtimes G) \otimes_{A \rtimes G_x} \big( {\ts \bigoplus_{g \in G_x} \C v_{gx} \big) = 
\bigoplus_{g \in G} \C v_{gx} = I_x \,.}
\]
c) The group $G$ permutes the $A$-weight spaces of $I_x$ transitively, and the $x$-weight 
space is just $\C [G_x]$. Hence the functor $\mr{Ind}_{A \rtimes G_x}^{A \rtimes G}$ 
provides a bijection between $G_x$-subrepresentations of  $\C [G_x]$ and 
$A \rtimes G$-submodules of $I_x$. Since $\C [G_x]$ is completely reducible, so is $I_x$.\\
d) Let $(\pi ,V)$ be an irreducible $A \rtimes G$-module. There exists an $x$ for which $V$ 
contains a nonzero $A$-weight vector $v_x$. Then
\begin{align*}
& \phi : I_x \to V , \\
& \phi ((a \otimes g) v_x) = \pi (a \otimes g) v_x
\end{align*}
is a surjective homomorphism of $A \rtimes G$-modules and hence $V \cong I_x / \ker \phi$ 
is a quotient of $I_x$. But $I_x$ is completely reducible, so $V$ is also isomorphic to a 
direct summand of $I_x$.\\
e) The map $\mr{Ind}_{A \rtimes G_x}^{A \rtimes G}$ above remains bijective on the level 
of isomorphism classes of irreducible modules. As is well known, the number of inequivalent 
irreducible constituents of the left regular representation $\C [G_x]$ equals the number of 
conjugacy classes in $G_x . \qquad \Box$
\\[2mm]

Extended quotients are very useful to describe the homology of crossed products, so us
recall these. Write
\begin{equation}\label{eq:1.1}
\tilde X = \{ (g,x) \in G \times X : g (x) = x \} .
\end{equation}
Then $G$ acts on $\tilde X$ by $g (g',x) = (g g' g^{-1},g(x))$, and the extended quotient
of $X$ by $G$ is defined as $\tilde X /G$. Write $X^g = \{ x \in X : g(x) = x \}$, let $Z_G (g)$ 
be the centralizer of $g$ in $G$ and let $G/ \sim$ be the collection of conjugacy classes in $G$. 
We can also construct the extended quotient as a disjoint union:
\begin{equation}\label{eq:1.2}
\begin{split}
\tilde X /G & = \Big( \bigcup_{g \in G} (g,X^g) \Big) / G = \bigcup_{c \in G / \sim} 
\Big( \bigcup_{g \in c} (g,X^g) / G \Big) \\
& \cong \bigcup_{g/ \sim \in G / \sim} (g,X^g / Z_G (g)) \cong 
\bigsqcup_{g/\sim \in G / \sim} X^g / Z_G (g) .
\end{split}
\end{equation}
Suppose now that $X$ is a nonsingular affine complex variety, and that $G$ acts on it by
algebraic isomorphisms. Then $X^g$ is an algebraic variety as well, and with the slice theorem
one deduces that $X^g$ is a smooth manifold. Hence $X^g$ and $\tilde X$ are nonsingular
affine varieties, although usually the latter has components of different dimensions.

For $A$ we take the algebra $\mc O (X)$ of regular functions on $X$.
Let $\Omega^n (X)$ denote the space of algebraic $n$-forms on $X$, and $H_{DR}^n (X)$ the
(algebraic) De Rham-cohomology of $X$.  The following result generalizes the
classical Hochschild--Kostant--Rosenberg theorem:

\begin{thm}\label{thm:1.2}
Let $X$ and $G$ be as above.
The Hochschild, cyclic and periodic cyclic homologies of $\mc O (X) \rtimes G$ are given by
\[
\begin{array}{ccl}
HH_n (\mc O (X) \rtimes G) & = & \Omega^n \big( \tilde X \big)^G , \\
HC_n (\mc O (X) \rtimes G) & = & 
\big( \Omega^n \big( \tilde X \big) / \textup{d} \Omega^{n-1} \big( \tilde X \big) \oplus 
H_{DR}^{n-2}\big( \tilde X \big) \oplus H_{DR}^{n-4}\big( \tilde X \big) \oplus \cdots \big)^G , \\
HP_n (\mc O (X) \rtimes G) & = & 
\big( {\ts \bigoplus_{m \in \Z}} H_{DR}^{n+2m} \big( \tilde X \big) \big)^G .
\end{array}
\]
\end{thm}
\emph{Proof.}
This result and various generalizations were proven by Brylinski and Nistor, 
see \cite{Bry} and \cite[Theorem 2.11]{Nis}. $\qquad \Box$
\\[2mm]

To make the above isomorphisms more explicit we need to recall the definitions of 
Hochschild and cyclic homology. For any unital algebra $A$ and any $A$-bimodule $M$ there 
is a differential complex $(C_* (A,M),b)$, where $C_n (A,M) = M \otimes A^{\otimes n}$ 
and $b : C_n (A,M) \to C_{n-1} (A,M)$ is defined on elementary tensors by
\begin{equation}\label{eq:1.3}
\begin{split}
b ( m \otimes a_1 \otimes \ldots \otimes a_n ) = \: & {\ts \sum_{i=1}^{n-1}} (-1 )^i 
m \otimes a_1 \otimes \ldots \otimes a_i a_{i+1} \otimes \ldots \otimes a_n \, + \\
& m a_1 \otimes a_1 \otimes \ldots \otimes a_n  + 
(-1 )^n a_n m \otimes a_1 \otimes \ldots \otimes a_{n-1} .
\end{split}
\end{equation}
This is known as the Hochschild complex, and its homology is denoted by $H_* (A,M)$.
When $M = A$, we omit it from the notation and get $H_n (C_* (A),b) = HH_n (A)$.

We can extend $C_* (A)$ to a so-called mixed complex $\mc B_* (A)$:
\begin{equation}\label{eq:mixedcomplex}
\begin{array}{ccccc}
\downarrow & & \downarrow & & \downarrow \\
A^{\otimes 3} & \xleftarrow{B} & A^{\otimes 2} & \xleftarrow{B} & A \\
\downarrow b & & \downarrow b & & \\
A^{\otimes 2} & \xleftarrow{B} & A \; \, \ & & \\
\downarrow b & & & & \\
A \; \, \ & & & & 
\end{array}
\end{equation}
The degree on this complex is the sum of the horizontal and vertical positions. The operator 
$B : A^{\otimes (n+1)} \to A^{\otimes (n+2)}$ satisfies $b B + B b = 0$ and is given explicitly by
\begin{equation}\label{eq:1.4}
B(a_1 \otimes \cdots \otimes a_n ) = \sum_{i=0}^n (-1 )^{ni} (1 \otimes a_i - a_i \otimes 1) 
\otimes a_{i+1} \otimes \cdots \otimes a_n \otimes a_1 \otimes \cdots \otimes a_{i-1} .
\end{equation}
The cyclic homology of $A$ is $HC_n (A) = H_n (\mc B_* (A),b,B)$.

The remainder of this section, which will only be used in Theorem \ref{thm:4.4} and not for 
Sections 4--6, provides a more detailed analysis of Theorem \ref{thm:1.2}.
For $A$ and $G$ as in \eqref{eq:1.5}, $C_n (A \rtimes G)$ has a subspace 
\[
C_n (A)_g := \mr{span} \{ g a_0 \otimes a_1 \otimes \cdots \otimes a_n : a_i \in A \} .
\]
Notice that $(C_* (A)_g ,b)$ is a subcomplex of $C_* (A \rtimes G) ,b)$. Its homology is
\[
H_n (C_* (A)_g ,b) = H_n (A ,A_g ) ,
\]
where the $A$-bimodule structure on $M = A_g$ is given by $a \cdot m \cdot a' = g^{-1}(a) m a'$.
According to \cite[Proposition 4.6]{GeJo} 
\begin{equation}\label{eq:1.6}
HH_n (A \rtimes G) \cong H_n \big( {\ts \bigoplus_{g \in G} C_* (A)_g ,b \big)_G \cong 
\big( \bigoplus_{g \in G}} H_n (A,A_g ) \big)^G ,
\end{equation}
where the sub- and superscripts $G$ mean coinvariants and invariants, respectively. The $G$-action
is defined via the inclusion $C_* (A)_g \to C_* (A \rtimes G)$:
\[
h \cdot (g a_0 \otimes a_1 \otimes \cdots \otimes a_n ) = 
(hgh^{-1}) \, h(a_0 ) \otimes h(a_1 ) \otimes \cdots \otimes h(a_n ) .
\]
Moreover for $A = \mc O (X)$ and $G$ as in Theorem \ref{thm:1.2}, \cite[Corollary 2.12]{Nis} says that
the inclusion $X^g \to X$ induces an isomorphism
\begin{equation}\label{eq:1.7}
H_n (\mc O (X) ,\mc O (X)_g ) \cong HH_n (\mc O (X^g )) \cong \Omega^n (X^g ) .
\end{equation}
The antisymmetrization map $\epsilon_n$ from \cite[Section 1.3]{Lod} gives a commutative diagram
\begin{equation}\label{eq:1.8}
\begin{array}{ccc}
\Omega^n (X^g ) & \longleftarrow & \Omega^n (\mc O (X)^g ) \\
\downarrow \epsilon_n & & \downarrow \epsilon_n \\
HH_n (\mc O (X^g )) & \longleftarrow & HH_n (\mc O (X )^g ) ,
\end{array}
\end{equation}
where $\Omega^n (B)$ denotes the module of K\"ahler $n$-forms over a commutative algebra $B$.
Since $X^g$ is nonsingular, the famous Hochschild--Kostant--Rosenberg theorem assures that the 
left $\epsilon_n$ is a bijection.

Since $X^g$ is closed in $X$, the restriction map $\mc O (X) \to \mc O (X^g )$ is surjective, and 
since $g$ has finite order it admits a linear splitting $s : \mc O (X^g ) \to \mc O (X)^g$. This
induces a splitting 
\[
\Omega^n (s) : \Omega^n (X^g ) \to \Omega^n (\mc O (X)^g ) 
\]
of the upper line of \eqref{eq:1.8}.
Now $\epsilon_n \circ \Omega^n (s) \circ \epsilon_n^{-1}$ in \eqref{eq:1.8} defines a splitting of
\begin{equation}\label{eq:1.9}
HH_n (\mc O (X )^g ) \to HH_n (\mc O (X^g )) ,
\end{equation}
which therefore is surjective. 
We can embed $C_* (\mc O (X)^g )$ in $C_* (\mc O (X))_g$, simply by writing a $g$ on the left.
Together with \eqref{eq:1.6}, \eqref{eq:1.7} and \eqref{eq:1.9} this shows that the inclusion 
\begin{equation}\label{eq:1.10} 
{\ts \bigoplus_{g \in G}} \, g C_* (\mc O (X)^g ) \to C_* (\mc O(X) \rtimes G)
\end{equation}
induces a surjection on homology. Moreover, because $HH_* (\mc O (X) \rtimes G)$ consists of 
$G$-coinvariants, it suffices to restrict the direct sum to one element in every conjugacy
class of $G$.

Let $\langle g \rangle \subset G$ be the cyclic group generated by $g \in G$, so that $g$ and 
$\mc O (X)^g$ lie in the commutative algebra 
$\mc O^g := \C [\langle g \rangle] \otimes \mc O (X)^g$. There is a natural surjection
\begin{align*}
& \pi_n : C_n (\mc O^g ) \to \Omega^n (\mc O^g ) , \\
& \pi_n (a_0 \otimes a_1 \otimes \cdots \otimes a_n ) = a_0 \, d a_1 \cdots d a_n .
\end{align*}
According to \cite[Lemma 1.3.14 and Proposition 2.3.4]{Lod} we have
\[
\pi_n \circ b = 0 \quad \mr{and} \quad \pi_{n+1} \circ B = (n+1) \, d \circ \pi_n .
\]
Therefore $\pi_n / n!$ induces a map from $(\mc B_* (\mc O^g ) ,b,B)$ to the mixed complex
\begin{equation}\label{eq:1.11}
\begin{array}{ccccc}
\downarrow & & \downarrow & & \downarrow \\
\Omega^2 (\mc O^g ) & \xleftarrow{d} & \Omega^1 (\mc O^g ) & 
\xleftarrow{d} & \Omega^0 (\mc O^g ) \\
\downarrow 0 & & \downarrow 0 & & \\
\Omega^1 (\mc O^g ) & \xleftarrow{d} & \Omega^0 (\mc O^g ) & & \\
\downarrow 0 & & & & \\
\Omega^0 (\mc O^g ) & & & & 
\end{array}
\end{equation}
Since $\Omega^n (\mc O (X)^g ) \to \Omega^n (\mc O (X^g )) = \Omega^n (X^g )$ is surjective, the 
comparison of Theorem \ref{thm:1.2} with the homology of \eqref{eq:1.11} shows that this gives us
a surjection
\[
HC_n (\mc O^g ) \to HC_n (\mc O (X^g )) .
\]
As $\C [\langle g \rangle]$ barely plays a role here, the map remains surjective when we restrict 
it to elements coming from $g \mc B_* (\mc O (X)^g ) \subset \mc B_* (\mc O^g )$. From this, 
Theorem \ref{thm:1.2} and \eqref{eq:1.10} we deduce the result that we were after:

\begin{lem}\label{lem:1.3} 
Let $X$ and $G$ be as in Theorem \ref{thm:1.2}. All elements of $HH_n (\mc O (X) \rtimes G)$ and
of $HC_n (\mc O (X) \rtimes G)$ can be represented by cycles in
$\bigoplus_{g \in G} g C_* (\mc O (X)^g )$. It also suffices to take the direct sum over
representatives of the conjugacy classes in $G$.
\end{lem}
\vspace{4mm}

\section{Graded Hecke algebras}
\label{sec:3}

For the construction of graded Hecke algebras we will use the following objects:
\begin{itemize}
\item a finite dimensional real inner product space $\mf a$,
\item the linear dual $\mf a^*$ of $\mf a$,
\item a (reduced) root system $R$ in $\mf a^*$,
\item the dual root system $R^\vee$ in $\mf a$,
\item a basis $\Pi$ of $R$.
\end{itemize}
We call 
\begin{equation}
\tilde{\mc R} = (\mf a^* ,R, \mf a, R^\vee, \Pi )
\end{equation}
a degenerate root datum. We neither assume that $R$ is crystallographic, nor 
that $\Pi$ spans $\mf a^*$. In fact $R$ is even allowed to be empty. 
Our degenerate root datum gives rise to
\begin{itemize}
\item the complexifications $\mf t$ and $\mf t^*$ of $\mf a$ and $\mf a^*$,
\item the symmetric algebra $S (\mf t^*)$ of $\mf t^*$,
\item the Weyl group $W$ of $R$,
\item the set $S = \{ s_\alpha : \alpha \in \Pi \}$ of simple reflections in $W$,
\item the complex group algebra $\mh C [W]$.
\end{itemize}
Choose formal parameters $\mb k_\alpha$ for $\alpha \in \Pi$, with the property that 
$\mb k_\alpha = \mb k_\beta$ if $\alpha$ and $\beta$ are conjugate under $W$. 
The graded Hecke algebra $\tilde{\mh H} (\tilde{\mc R})$ corresponding to 
$\tilde{\mc R}$ is defined as follows. As a complex vector space
\[
\tilde{\mh H} (\tilde{\mc R}) = \mh C [W] \otimes S (\mf t^* ) \otimes 
\mh C [\{ \mb k_\alpha : \alpha \in \Pi \} ] .
\]
The multiplication in $\tilde{\mh H} (\tilde{\mc R})$ is determined by the following rules:
\begin{itemize}
\item $\mh C[W] \,, S (\mf t^* )$ and $\mh C[\{ \mb k_\alpha : \alpha \in \Pi \} ]$ 
are canonically embedded as subalgebras,
\item the $\mb k_\alpha$ are central in $\tilde{\mh H} (\tilde{\mc R})$,
\item for $x \in \mf t^*$ and $s_\alpha \in S$ we have the cross relation
\begin{equation}\label{eq:3.1}
x \cdot s_\alpha - s_\alpha \cdot s_\alpha (x) = 
\mb k_\alpha \inp{x}{\alpha^\vee} \,.
\end{equation}
\end{itemize}
We define a grading on $\tilde{\mh H} (\tilde{\mc R})$ by requiring that $\mf t^*$
and the $\mb k_\alpha$ are in degree one, while $W$ has degree zero.

In fact we will only study specializations of this algebra.
Pick complex numbers $k_\alpha \in \mh C$ for $\alpha \in \Pi$, such 
that $k_\alpha = k_\beta$ if $\alpha$ and $\beta$ are conjugate under
$W$. Let $\mh C_k$ be the onedimensional $\mh C[\{ \mb k_\alpha : \alpha \in \Pi \} ]$-module 
on which $\mb k_\alpha$ acts as multiplication by $k_\alpha$. We define 
\begin{equation}
\mh H = \mh H (\tilde{\mc R},k) = \tilde{\mh H} (\tilde{\mc R})
\otimes_{\mh C[\{ \mb k_\alpha : \alpha \in \Pi \} ]} \mh C_k
\end{equation}
With some abuse of terminology $\mh H (\tilde{\mc R},k)$ is also
called a graded Hecke algebra. Notice that as a vector space
$\mh H (\tilde{\mc R},k)$ equals $\mh C[W] \otimes S (\mf t^*)$, and
that the cross relation \eqref{eq:3.1} now holds with $\mb k_\alpha$
replaced by $k_\alpha$:
\begin{equation}\label{eq:3.2}
x \cdot s_\alpha - s_\alpha \cdot s_\alpha (x) = k_\alpha \inp{x}{\alpha^\vee} .
\end{equation}
More generally, for any $s_\alpha \in S$ and $p \in S (\mf t^* )$ we have
\begin{equation}\label{eq:3.5}
p \cdot s_\alpha - s_\alpha \cdot s_\alpha (p) = k_\alpha \frac{p - s_\alpha (p)}{\alpha} .
\end{equation}
An automorphism $\gamma$ of the Dynkin diagram of the based root system $(R,\Pi )$ 
is a bijection $\Pi \to \Pi$ such that
\begin{equation}
\inp{\gamma (\alpha )}{\gamma (\beta )^\vee} = \inp{\alpha}{\beta^\vee} 
\qquad \forall \alpha ,\beta \in \Pi \,.
\end{equation}
Such a $\gamma$ naturally induces automorphisms of $R, R^\vee$ and $W$.
Moreover we will assume that the action of $\gamma$ on the linear span of the coroots 
is extended somehow to an orthogonal endomorphism of $\mf t$. If $\gamma$ and $\gamma'$
act in the same way on $R$ but differently on $\mf t$, then we will sloppily regard 
them as different diagram automorphisms.

Let $\Gamma$ be a finite group of diagram automorphisms of $(R,\Pi)$. Groups like
\begin{equation}
W' := \Gamma \ltimes W
\end{equation}
typically arise from larger Weyl groups as the isotropy groups of points in some torus, or as 
normalizers of some parabolic subgroup \cite{How}.
Assume that $k_{\gamma (\alpha )} = k_\alpha \; \forall \alpha \in \Pi , \gamma \in \Gamma$.
Then $\Gamma$ acts on $\mh H$ by the algebra homomorphisms
\begin{equation}\label{eq:3.6}
\begin{split}
& \psi_\gamma : \mh H \to \mh H \,, \\
& \psi_\gamma (x s_\alpha ) = \gamma (x) s_{\gamma (\alpha )} \qquad x \in \mf t^* , \alpha \in \Pi \,.
\end{split}
\end{equation}
Thus we can form the crossed product
\begin{equation}\label{eq:3.4}
\mh H' := \Gamma \ltimes \mh H = \Gamma \ltimes \mh H (\tilde{\mc R},k ) \,,
\end{equation}
which we call an extended graded Hecke algebra.

We define a $\Z$-grading on $\mh H'$ by deg$(w) = 0 \; \forall w \in W'$ and
deg$(x) = 1 \; \forall x \in \mf t^*$. \linebreak However, the algebra $\mh H'$ is in
general not graded, only filtered. That is, the product $h_1 h_2$ of two homogeneous elements 
$h_1 ,h_2 \in \mh H'$ need not be homogeneous, but all its homogeneous components
have degree at most $\text{deg}(h_1 ) + \text{deg}(h_2 )$. More precisely, from \eqref{eq:3.5} 
we see that the part of $h_1 h_2$ that depends on the parameters $k_\alpha$ has degree strictly 
lower than $\text{deg}(h_1 ) + \text{deg}(h_2 )$.

Let us mention some special cases in which 
$\mh H' = \Gamma \ltimes \mh H (\tilde{\mc R},k)$ is graded:
\begin{itemize}
\item if $R = \es$ then $\mh H' = \Gamma \ltimes \tilde{\mh H} (\tilde{\mc R}) = 
\Gamma \ltimes S (\mf t^* )$,
\item if $k_\alpha = 0 \; \forall \alpha \in \Pi$ then $\Gamma \ltimes \mh H (\tilde{\mc R},k)$ 
is the crossed product $W' \ltimes S (\mf t^* )$, with the cross relations
\begin{equation}
w \cdot x = w(x) \cdot w \qquad w \in W' , x \in \mf t^* .
\end{equation}
\end{itemize}
Multiplication with any $z \in \mh C^\times$ defines a bijection $m_z : \mf t^* \to \mf t^*$,
which clearly extends to an algebra automorphism of $S(\mf t^* )$. From the cross relation
\eqref{eq:3.2} we see that it extends even further, to an algebra isomorphism
\begin{equation}\label{eq:3.3}
m_z : \Gamma \ltimes \mh H (\tilde{\mc R},zk) \to \Gamma \ltimes \mh H (\tilde{\mc R}, k)
\end{equation}
which is the identity on $\mh C[W']$. Notice that the homomorphism \eqref{eq:3.3} remains 
well-defined for $z=0$, only then it ceases to be bijective.

In particular, if all $\alpha \in R$ are conjugate under $W'$, then there are essentially 
only two graded Hecke algebras attached to $(\tilde{\mc R},\Gamma)$: one with $k=0$ and 
one with $k \neq 0$.
\vspace{4mm}

\section{The homology of graded Hecke algebras}
\label{sec:4}

In this section we will prove that the periodic cyclic homology of a graded Hecke algebra
is the same as that of its underlying finite Weyl group. In fact we immediately generalize
this to the extended graded Hecke algebras from \eqref{eq:3.4}, since this does not make
the proof more difficult. Throughout we will use the notations from Section \ref{sec:3}.

\begin{lem}\label{lem:4.1}
Consider the extended graded Hecke algebra 
$\Gamma \ltimes \mh H (\tilde{\mc R}, 0) = W' \ltimes S (\mf t^* )$ with parameter $k=0$.
\begin{description}
\item[a)] $HH_n \big( W' \ltimes S (\mf t^* ) \big) = 0$ for all $n > \dim_\C (\mf t^* )$.
\item[b)] $HC_n \big( W' \ltimes S (\mf t^* ) \big) = 
HP_n \big( W' \ltimes S (\mf t^* ) \big) = 0$ for all odd $n > \dim_\C (\mf t^* )$.
\item[c)] The inclusion $\C [W'] \to W' \ltimes S (\mf t^* )$ induces an isomorphism
on \pch .
\end{description}
\end{lem}
\emph{Proof.}
a) follows immediately from Theorem \ref{thm:1.2}.\\
b) Every component of $\widetilde{\mf t}$ is a vector space, and in particular 
contractible as a topological space. Therefore $H_{DR}^n \big( \widetilde{\mf t} \big) = 0$
for all $n > 0$. Now apply Theorem \ref{thm:1.2}.\\
c) By the above we only have to show this in degree 0. Again by Theorem \ref{thm:1.2}
\[
\begin{split}
HP_0 \big( W' \ltimes S (\mf t^* ) \big) & \cong H_{DR}^0 \big( \widetilde{\mf t} \big)^{W'} 
\cong \big( {\ts \bigoplus_{w \in W'}} H_{DR}^0 \big( \mf t^w \big) \big)^{W'} \\
& \cong \big( \C [W'] \big)^{W'} \cong HP_0 \big( \C [W'] \big)^{W'} .
\end{split}
\]
According to \cite[Theorem 2.11]{Nis} all the above isomorphisms can be made natural. In
particular $\C [W'] \to W' \ltimes S (\mf t^* )$ induces an isomorphism on $HP_* . \qquad \Box$
\\[2mm]

To transfer these results to $\mh H' = \Gamma \ltimes \mh H (\tilde{\mc R},k)$ with 
$k \neq 0$, we will use spectral sequences. For $d \in \Z_{\geq 0}$ let $\mh H'_{\leq d}$ 
denote the vector space of elements of degree $\leq d$ in $\mh H'$. 
We generalize this filtration to the tensor power $(\mh H' )^{\otimes n}$ by
\begin{equation}\label{eq:4.1}
(\mh H')^{\otimes n}_{\leq d} := \sum_{d_1 + \cdots + d_n \leq d} \mh H'_{\leq d_1} \otimes
\cdots \otimes \mh H'_{\leq d_n} .
\end{equation}
Equivalently one can take the sum over all $n$-tuples $(d_1 ,\ldots ,d_n )$ of nonnegative 
integers with $d_1 + \cdots + d_n = d$.

\begin{lem}\label{lem:4.2}
$HH_n \big( \Gamma \ltimes \mh H (\tilde{\mc R},k) \big) = 0 \quad \mr{for} \; n > \dim_\C (\mf t^* )$.
\end{lem}
\emph{Proof.}
The Hochschild complex $(C_* (\mh H') ,b)$ from \eqref{eq:1.3} is filtered by
\eqref{eq:4.1}. That is, we put 
\[
F_p C_n (\mh H' )= ( \mh H' )^{\otimes (n+1)}_{\leq p} .
\]
This gives rise to a spectral sequence converging to
$HH_* (\mh H' )$, see for example \cite[Chapter XV]{CaEi}. Its first term is 
\[
E^1_{p,q} = H_{p+q} \big( F_p C_* (\mh H') / F_{p-1} C_* (\mh H') \big) . 
\]
The boundary map on $F_p C_* (\mh H')/ F_{p-1} C_* (\mh H')$ is given by \eqref{eq:1.3}, 
but we may ignore all terms that are not in top degree. From equations \eqref{eq:3.2} and 
\eqref{eq:3.6} we see that the resulting map is actually independent of $k$. 
So we can determine $E^1_{p,q}$ for all parameters, just by looking at the case $k = 0$.

But for $k=0$ the algebra $\Gamma \ltimes \mh H (\tilde{\mc R},0) = W' \ltimes S(\mf t^* )$ is
really graded, so filtering the Hochschild complex does not add anything, and the spectral
sequence stabilizes already at $E^1_{*,*}$. From Lemma \ref{lem:4.1}.a we deduce that 
$E^1_{p,q} = 0$ if $p+q > \dim_\C (\mf t^* )$.

For general parameters $k_\alpha$ we cannot say immediately whether our spectral sequence 
stabilizes at the first term, but in any case $E^r_{p,q}$ is a subquotient of $E^1_{p,q}$. 
Hence $E^\infty_{p,q} = 0$ whenever $p+q > \dim_\C (\mf t^* )$, and the lemma follows.
$\qquad \Box$
\vspace{3mm}

As a corollary of Lemma \ref{lem:4.2} we find that for $n > \dim_\C (\mf t^* )$:
\begin{equation}\label{eq:4.3}
HC_n \big( \Gamma \ltimes \mh H (\tilde{\mc R},k) \big) \cong 
HP_n \big( \Gamma \ltimes \mh H (\tilde{\mc R},k) \big) .
\end{equation}
Indeed, this can be seen from Connes' periodicity exact sequence, which relates Hochschild and
cyclic homology. On the other hand, we can also try to compute the cyclic homology with a
spectral sequence, and this leads to the following result.

\begin{thm}\label{thm:4.3} 
For $n > \dim_\C (\mf t^* )$ we have
\[
HP_n \big( \Gamma \ltimes \mh H (\tilde{\mc R},k) \big) \cong 
HC_n \big( \Gamma \ltimes \mh H (\tilde{\mc R},k) \big) \cong HC_n \big( \C [\Gamma \ltimes W] \big) .
\]
The inclusion $\C [\Gamma \ltimes W] \to \Gamma \ltimes \mh H (\tilde{\mc R},k)$ induces an
isomorphism on periodic cyclic homology.
\end{thm}
\emph{Proof.}
Recall from \eqref{eq:mixedcomplex} that $HC_n (\mh H')$ is computed by the mixed complex
$(\mc B_* (\mh H'),b,B)$. Its space of $n$-chains is
\[
\mc B_n (\mh H' ) := \mh H'^{\otimes (n+1)} \oplus \mh H'^{\otimes n} \oplus \cdots \oplus \mh H' .
\]
We introduce a filtration on this differential complex by
\[
F_p \mc B_n (\mh H' ) := ( \mh H' )^{\otimes (n+1)}_{\leq p} \oplus (\mh H' )^{\otimes n}_{\leq p} 
\oplus \cdots \oplus \mh H'_{\leq p} .
\]
As in the previous proof, this yields a spectral sequence $E^r_{p,q}$ converging to
$HC_{p+q}(\mh H' )$. In particular
\begin{equation}\label{eq:4.6}
E^1_{p,q} = H_{p+q} \big( F_p \mc B_* (\mh H' ) / F_{p-1} \mc B_* (\mh H' ) \big) ,
\end{equation}
By \eqref{eq:3.2} and \eqref{eq:3.6} the boundary maps in 
$F_p \mc B_* (\mh H' ) / F_{p-1} \mc B_* (\mh H' )$ are independent of the parameters $k_\alpha$. 
Hence the vector spaces $E^1_{p,q}$ do not depend on $k$, and we can determine them from the 
special case $k=0$.

Because the algebra $\Gamma \ltimes \mh H (\tilde{\mc R},0) = W' \ltimes S(\mf t^* )$ is
graded, the spectral sequence for this algebra already stabilizes at $E^1_{*,*}$, and
$E^1_{p,q}$ is just the degree $p$ part of $HC_{p+q} (W' \ltimes S(\mf t^* ) )$.
Unlike Hochschild homology, $HC_n (W' \ltimes S(\mf t^* ) )$ does not vanish for (sufficiently)
large $n$, but from Lemma \ref{lem:4.1}.b we do know that $HC_n (W' \ltimes S(\mf t^* ) ) = 0$
when $n > \dim_\C (\mf t^* )$ is odd. Moreover, for even $n > \dim_\C (\mf t^* )$ 
Theorem \ref{thm:1.2} says that
\[
HC_n \big( W' \ltimes S(\mf t^* ) \big) = HP_n \big( W' \ltimes S(\mf t^* ) \big) .
\]
It follows from Lemma \ref{lem:4.1}.c that, for all $n > \dim_\C (\mf t^* ) ,\, 
HC_n (W' \ltimes S(\mf t^* ) )$ has no parts in degrees $p>0$, and that its degree $p = 0$ 
part is $HH_0 \big( \C [W'] \big)$.

Now we return to general $k$ and consider various $p$ and $q$. By definition $E^r_{p,q} = 0$ 
if $p < 0$ or $p+q < 0$. We can not yet say much about $E^1_{p,q}$ when 
$0 \leq p+q \leq \dim_\C (\mf t^* )$, but we can do without. 

Finally, we pick $p,q \in \Z$ such that $p+q > \dim_\C (\mf t^* )$. By the above $E^1_{p,q} = 0$ 
unless $p=0$ and $q$ is even, in which case $E^1_{p,q} = HH_0 \big( \C [W'] \big)$. For every 
$r \in \Z_{\geq 0}$ there is a boundary map $\partial^r_{p,q} : E^r_{p,q} \to E^r_{p-r,q+r-1}$,
and $E^{r+1}_{*,*}$ is the homology of $\big( E^r_{*,*} , \partial^r_{*,*} \big)$. We claim that
$\partial^r_{p,q} = 0$ whenever $r \geq 1$ and $p+q > \dim_\C (\mf t^* )$. Indeed, for $p > 0$ 
the domain is zero, while for $p=0$ the range is $E^r_{-r,q+r-1}$, which is zero because $-r < 0$.

We conclude that in the range $p+q > \dim_\C (\mf t^* )$ our spectral sequence stabilizes at
$r=1$, for all $k$. Moreover, we already showed that the vector spaces $E^1_{p,q}$ do not depend
on $k$. Taking the limit $r \to \infty$ and using the convergence, we find that 
$HC_{p+q} (\mh H' )$ does not depend on $k$. In view of \eqref{eq:4.3} the same goes for
$HP_{p+q} (\mh H' )$. But the functor $HP_*$ is 2-periodic, so actually $HP_n (\mh H' )$ is 
independent of $k$ for all $n \in \Z$.
Together with Lemma \ref{lem:4.1}.c this shows that 
\[
HP_n \big( \Gamma \ltimes \mh H (\tilde{\mc R},k) \big) \cong 
HP_n \big( \C [\Gamma \ltimes W] \big) \qquad \forall n \in \Z .
\]
Furthermore we noticed that for large $n$ the representatives of 
\[
HC_n \big( \Gamma \ltimes \mh H (\tilde{\mc R},k) \big) \cong
HP_n \big( \Gamma \ltimes \mh H (\tilde{\mc R},k) \big)
\]
all lie in the image of 
\begin{equation}\label{eq:4.5}
HP_n \big( \C [\Gamma \ltimes W] \big) \to HP_n \big(\Gamma \ltimes \mh H (\tilde{\mc R},k) \big) .
\end{equation}
Therefore \eqref{eq:4.5} is a linear bijection. $\qquad \Box$
\\[2mm]

Our next result is a clear improvement on Lemma \ref{lem:4.2} and Theorem \ref{thm:4.3}. 
Although the proofs of these three results could have been combined in one, we decided
against this. The previous results are relatively general and the proofs probably can be
applied to similar algebras as well. The upcoming result however uses more subtle properties
of root systems and crossed products.

\begin{thm}\label{thm:4.4}
For all $n \in \Z_{\geq 0}$ there are isomorphisms
\[
\begin{array}{ccc}
HH_n (\Gamma \ltimes \mh H (\tilde{\mc R},k)) & \cong & HH_n (W' \ltimes S(\mf t^* )) , \\
HC_n (\Gamma \ltimes \mh H (\tilde{\mc R},k)) & \cong & HC_n (W' \ltimes S(\mf t^* )) , \\
HP_n (\Gamma \ltimes \mh H (\tilde{\mc R},k)) & \cong & HP_n (W' \ltimes S(\mf t^* )) . \\
\end{array}
\]
\end{thm}
\emph{Proof.}
In the proofs of Lemma \ref{lem:4.2} and Theorem \ref{thm:4.3} we constructed spectral
sequences $E^r_{p,q}$ converging to $HH_{p+q}(\Gamma \ltimes \mh H (\tilde{\mc R},k))$ and
$HC_{p+q}(\Gamma \ltimes \mh H (\tilde{\mc R},k))$, respectively. We showed that the spaces
$E^1_{p,q}$ do not depend on $k$, but it is conceivable that boundary maps $\partial^r_{p,q}$ 
do. With a clever choice of representatives for the homology classes we will show that
$\partial^r_{p,q} = 0$ for all $r \geq 1$. We will only write this down for $HC_*$, the 
Hochschild homology can be handled in the same way. Recall from page \pageref{eq:4.6} that
\begin{equation}\label{eq:4.7}
{\ts \bigoplus_{p+q = n}} E^1_{p,q} \cong HC_n (W' \ltimes S (\mf t^* )) .
\end{equation}
By Lemma \ref{lem:1.3} every element of $HC_n (W' \ltimes S (\mf t^* ))$ can be represented
by a cycle in 
\[
\bigoplus_{g \in W'} g C_* \big( \mc O (\mf t )^g \big) = 
\bigoplus_{g \in W'} g \bigoplus_{m \geq 1} \big( \mc O (\mf t )^g \big)^{\otimes m} .
\]
Moreover we need only one $g$ from every conjugacy class in $W'$, and the map from cycles
to homology classes factors through $\mc O (\mf t )^g \to \mc O (\mf t^g )$.

Recall that the closed Weyl chamber
\[
\mf a^+ := \{ \lambda \in \mf a : \inp{\alpha}{\lambda} \geq 0 \; \forall \alpha \in \Pi \}
\]
is a fundamental domain for the action of $W$ on $\mf a$. Take any $g \in W'$ and write
\begin{align*}
& R_g = \{ \alpha \in R : \alpha \perp \mf a^g \} , \\ 
& \mf a^{R_g} = \{ \lambda \in \mf a : \inp{\alpha}{\lambda} = 0 \; \forall \alpha \in R_g \} .
\end{align*}
Clearly $\mf a^g \subset \mf a^{R_g}$.
By \cite[Theorem 1.12]{Hum} there exists $u \in W$ such that 
\[
u (\mf a^{R_g} ) = \mf a^P := 
\{ \lambda \in \mf a : \inp{\alpha}{\lambda} = 0 \; \forall \alpha \in P \}
\]
for some set $P \subset \Pi$ of simple roots. Replacing $g$ by $u g u^{-1}$, we may assume
that $\mf a^g \subset \mf a^P$. Then
\[
\mf a^g = (\mf a^g \cap \mf a^+ ) + (\mf a^g \cap - \mf a^+ ).\
\]
Since $\Gamma$ permutes $\Pi$, it stabilizes $\mf a^+$. So if we write $g = \gamma w$ with
$\gamma \in \Gamma$ and $w \in W$, then 
\[
w (\mf a^g \cap \mf a^+ ) = \gamma^{-1} (\mf a^g \cap \mf a^+ ) \subset \mf a^+ .
\]
Hence $w$ must fix $\mf a^g \cap \mf a^+$, and we deduce that both $w$ and $\gamma$ fix
$\mf a^g$ pointwise. Moreover by \cite[Theorem 1.12]{Hum} $w$ lies in the parabolic 
subgroup $W_P \subset W$.
Let $\mf t_P \subset \mf t$ be the complex span of $\{ \alpha^\vee : \alpha \in P \}$. This
is a complement to $\mf t^P = \mf a^P + i \mf a^P$ in $\mf t$ and it is stable under $\gamma$
and $W_P$. Clearly the natural map $\mc O (\mf t / \mf t_P ) \to \mc O (\mf t^g )$ is 
surjective. Hence every element of $HC_n (W' \ltimes S (\mf t^* ))$ can be represented by a
cycle
\[
z \in {\ts \bigoplus_g} \, g C_* (\mc O (\mf t /\mf t_P )) ,
\]
where we sum only over $g = \gamma w$ as above. Now consider $z$ as an element of the mixed
complex $\mc B_* (\Gamma \ltimes \mh H (\tilde{\mc R},k))$. Since $\gamma$ stabilizes $\mf t_P$
and all the simple reflections $s_\alpha$ with $\alpha \in P$ act trivially on $\mf t / \mf t_P$,
the multiplication of $g = \gamma w \in \Gamma W_P$ with elements of $\mc O (\mf t / \mf t_P )$
does not depend on $k$. 
Consequently $z$ is also a cycle in $\mc B_* (\Gamma \ltimes \mh H (\tilde{\mc R},k))$ and
$\partial^r_{*,*}(z) = 0$ for all $r \geq 0$. 

By \eqref{eq:4.7} every element of $E^1_{*,*}$ can be obtained in this way, so 
$\partial^r_{*,*} = 0$ for all $r \geq 1$ and the spectral sequence $E^r_{*,*}$ stabilizes
at $r = 1$. This proves the theorem for $HH_*$ and $HC_*$, the statement for $HP_*$ was
already contained in Theorem \ref{thm:4.3}. $\qquad \Box$
\\[2mm]

We remark that 
\[
HH_n (\Gamma \ltimes \mh H (\tilde{\mc R},k)) \cong HH_n (W' \ltimes S(\mf t^* )) 
\]
is not an isomorphism of $S (\mf t^* )^{W'}$-modules. For $k \neq 0$ this structure is
distorted in the process of choosing suitable representatives of homology classes. 
However, it is still a little more than just a linear bijection between vector spaces.
The direct summand $\Omega^n (\mf t^g )^{Z_{W'}(g)}$ of $HH_n (W' \ltimes S (\mf t^* ))$
is a module over $\mc O (\mf t^g )^{Z_{W'}(g)} \subset S (\mf t^* )^{W'}$, and this module 
structure persists to the corresponding summand of $HH_n (\Gamma \ltimes \mh H (\tilde{\mc R},k))$.
\vspace{4mm}

\section{Periodic cyclic homology as sheaf cohomology}
\label{sec:2}

To draw representation theoretic consequences from knowledge of the \pch \: of an algebra,
we will relate it to its primitive ideal spectrum. This is possible for algebras of
finite type, as studied in \cite{KNS,BaNi}. 

Let $X$ be a complex, affine, algebraic variety
and $\mc O (X)$ its coordinate ring. A finite type $\mc O(X)$-algebra is an algebra $A$ 
together with a morphism from $\mc O (X)$ to the center of the multiplier algebra of $A$,
which makes $A$ into an $\mc O (X)$-module of finite rank. An algebra is said to be of 
finite type if it is a finite type $\mc O (X)$-algebra for some $X$. Notice that a unital
finite type algebra is always of finite type over its center.

Recall that a primitive ideal of an algebra $A$ is the annihilator of an irreducible $A$-module.
We endow the collection Prim$(A)$ of primitive ideals with the Jacobson topology, 
whose closed sets are of the form
\[
V(S) := \{ I \in \mr{Prim}(A) : S \subset I \} ,
\]
for any subset $S \subset A$. The resulting topological space is called the primitive ideal 
spectrum of $A$. Notice that Prim$(\mc O (X))$ is $X$ with the Zariski topology.

\begin{lem}\label{lem:PrimT1}
The primitive ideal spectrum of a finite type algebra is a $T_1$-space. 
\end{lem}
\emph{Proof.}
Recall that a topological space is $T_1$ if and only if every one point subset is closed.
Let $A$ be (not necessarily unital) finite type $\mc O (X)$-algebra, let $V$ be an irreducible 
$A$-module and $J \in \text{Prim}(A)$ the annihilator of $V$. 
We want to show that $V(J) = \{J\}$.

The $A$-module structure on $V$ extends naturally to the multiplier algebra $\mc M (A)$,
because $V = A V$. Let $\phi_A : \mc O (X) \to \mc M (A)$ be the morphism that defines the
$\mc O(X)$-algebra structure on $A$. It enables us to construct the unital finite 
$\mc O (X)$-algebra $A^+ := A \oplus \mc O (X)$ with multiplication 
\[
(a,f) (a' ,f') = (a a' + \phi_A (f) a' + a \, \phi_A (f') , f f' ) . 
\]
The $A$-module structure on $V$ extends to $A^+$ via 
\[
(a,f) v := (a + \phi_A (f)) v .
\]
Since $A^+$ has (at most) countable dimension over $\C$, the dimension of $V$ is countable 
or finite. Therefore Schur's lemma applies, and it tells us that $Z(A^+)$ acts on $V$ via a 
character $\chi : Z(A^+) \to \C$. Since $A^+$ is a unital finite type algebra, it has finite rank 
over $Z(A^+)$, and consequently the image of $A^+ \to \mr{End}_\C (V)$ has finite dimension. 
But $V$ is irreducible, so dim$(V)$ has to be finite. If $J^+$ is the annihilator of $V$ in $A^+$, 
then Wedderburn's theorem says that $A^+/J^+ \to \mr{End}_\C (V)$ is an isomorphism. 

Suppose that there is a primitive ideal $I$ of $A$ such that $I \supsetneq J$. 
Then $(I + J^+) / J^+$ is a nonzero ideal of the simple algebra $A^+/J^+ \cong \mr{End}_\C (V)$, 
so $(I + J^+) / J^+ = A^+ / J^+$. But then
\[
A = (I + J^+) \cap A = I + (J^+ \cap A) = I + J = I,
\]
which contradicts the primitivity of $I$. We conclude that the subset $V(J) \subset 
\mr{Prim}(A)$, which by definition is closed, equals $\{J\}. \qquad \Box$
\\[2mm]

For any unital finite algebra $A$ with center $\mc O (X)$ the map
\begin{equation}
\begin{split}
& \theta : \mr{Prim}(A) \to \mr{Prim}(\mc O (X)) \cong X ,\\
& \theta (I) = I \cap \mc O (X) 
\end{split}
\end{equation}
is a finite-to-one continuous surjection \cite[Lemma 1]{KNS}. In other words, Prim$(A)$ 
may be complicated, but cannot be too different from an algebraic variety. The \pch \:
of commutative finite type algebras was determined in \cite[Theorem 9]{KNS}:

\begin{thm}\label{thm:2.1}
Let $I \subset \mc O (X)$ be an ideal and $Y \subset X$ its zero locus. There is a natural
isomorphism
\[
HP_n (I) \cong {\ts \bigoplus_{m \in \Z}} \check H^{n+2m} (X^{an},Y^{an} ; \C ) .
\]
\end{thm}
Here the superscript ``an'' means that we endow this variety with the analytic topology,
and $\check H^*$ is (relative) \v Cech-cohomology.

It turns out that the \pch \ of a general finite type algebra depends only on its primitive
ideal spectrum. Following \cite{BaNi} we call a morphism $\phi : A \to B$ of finite type 
$\mc O (X)$-algebras spectrum preserving if:
\begin{itemize}
\item for each $J \in \mr{Prim}(B)$ there is exactly one $I \in \mr{Prim}(A)$ 
containing $\phi^{-1}(J)$,
\item the map Prim$(B) \to \mr{Prim}(A) : J \mapsto I$ is a bijection.
\end{itemize}

\begin{thm}\label{thm:2.2} \textup{\cite[Theorem 8]{BaNi}}\\
Let $\phi : A \to B$ be a spectrum preserving morphism of finite type $\mc O (X)$-algebras.
Then $HP_* (\phi ) : HP_* (A) \to HP_* (B)$ is an isomorphism.
\end{thm}

In view of these two Theorems, $HP_* (A)$ can be regarded as a kind of cohomology of Prim$(A)$.
But it requires some care to describe exactly how $HP_* (A)$ can be determined from the 
primitive ideal spectrum of $A$. Let Prim$(A)^{an}$ be the set Prim$(A)$ with the coarsest 
topology that makes
\[
\mr{Prim}(A)^{an} \to \mr{Prim}(A) \times X^{an} : I \mapsto (I, \theta (I))
\]
continuous. In this space every point has an open neighborhood homeomorphic to 
$\R^n \times \R_{\geq 0}^m$ for some $n,m \in \Z_{\geq 0}$, so roughly speaking it is a 
non-Hausdorff version of an orbifold.

If we would use a classical cohomology theory on Prim$(A)^{an}$, the homotopy axiom would imply
that we cannot see the entire space, only some Hausdorff quotient. The solution is to consider
not Prim$(A)^{an}$ but $X^{an}$, and to encode the additional structure of Prim$(A)$ in a sheaf
$S_A$ over $X^{an}$. This dictates that the stalk $S_A (x)$ must be the complex vector space
$\C \theta^{-1}(x)$, whose basis can be identified with the inequivalent irreducible $A$-modules 
with $\mc O(X)$-character $x \in X$. 

Sections of this sheaf will be modelled by idempotents in suitable completions of the algebra $A$.
More precisely, given an open $U \subset X^{an}$, pick a connected component $C$ of 
$\theta^{-1}(U) \subset \mr{Prim}(A)^{an}$, and let $C^H$ be its maximal Hausdorff quotient.
We note that $\theta : C \to U$ factors as 
\[
C \xrightarrow{q_C} C^H \xrightarrow{\theta_C} U .
\]
Let $s_C : C^H \to C$ be a section of $q_C$ and put
\begin{equation}\label{eq:2.2}
\begin{split}
& s : U \to {\ts \bigsqcup_{x \in U} } S_A (x) , \\
& s(x) = \sum_{c \in C^H : \theta_C (c) = x} s_C (c) .
\end{split}
\end{equation}
\begin{defn}\label{def:SA}
Let $A$ be a finite type $\mc O (X)$-algebra. Then $S_A$ is the sheaf over $X^{an}$ with stalks 
$S_A (x) = \C \theta^{-1} (x)$, whose continuous sections can be written locally as finite sums 
$\sum_i z_i s_i$, with $z_i \in \C$ and $s_i$ of the form \eqref{eq:2.2}.
\end{defn}

This construction is functorial in $A$. Indeed, suppose that $\phi : A \to B$ is a morphism of
finite type $\mc O (X)$-algebras. Then $V \mapsto B \otimes_A V$ is a functor from left 
$A$-modules to left $B$-modules, and it induces a map $S_A \to S_B$ of sheaves over $X$.
\\[2mm]
\textbf{Examples.}\\ 
Suppose that $A \subset \mc O (X) \otimes \mr{End}_\C (V)$ is unital and $e \in A$ 
is idempotent. Then we can get a global section $s_e \in S_A (X)$ as follows. 
For any $x \in X$ let 
\[
0 \subset V_x^1 \subset \cdots \subset V_x^d = V_x
\]
be a composition series of the $A$-module $V_x$, the vector space $V$ on which $\mc O(X)$ 
acts by $x$. For a primitive ideal $I \in \theta^{-1}(x)$ let 
\[
n(e,I) = {\ts \sum_m} \mr{rank} \big( e, V_x^m / V_x^{m-1} \big) \in \Z_{\geq 0} ,
\]
where the sum runs over all $m$ for which $I \subset A$ is the 
annihilator of the irreducible module $V_x^m / V_x^{m-1}$. Then 
\begin{equation}\label{eq:2.1}
s_e (x) := {\ts \sum_{I \in \theta^{-1} (x)} }  n(e,I) \, I
\end{equation}
defines a continuous section of $S_A$.

For a more concrete example, fix a point $x_0 \in X$ and consider the algebra
\[
B := \{ b \in \mc O (X) \otimes M_2 (\C) : b(x_0) \text{ is a diagonal matrix} \} .
\]
Then Prim$(B)$ is $X$, but with two points $x'_0 ,x''_0$ instead of $x_0$. The topology of 
Prim$(B)^{an}$ is such that every sequence $(x_n )_{n=1}^\infty$ in $X \setminus \{x_0\}$, 
which converges to $x_0$ in $X^{an}$, converges to both $x'_0$ and $x''_0$ in 
Prim$(B)^{an}$. The stalk $S_B (x)$ is canonically isomorphic to $\C$ for 
$x \in X \setminus \{x_0\}$, while $S_B (x_0) \cong \C^2$. In particular the sheaf $S_B$ on $X^{an}$
is not a vector bundle. A global section $s$ of $S_A$ is continuous if and only if the function
\[
X^{an} \to \C : x \mapsto \left\{ \begin{array}{lll}
s(x) & \text{if} & x \in X \setminus \{x_0\} \\
z_1 + z_2 & \text{if} & x = x_0 \;\; \text{and} \;\; s(x_0) = (z_1,z_2)
\end{array}
\right.
\]
is locally constant. Thus the continuous sections with $s(x_0) = (z_1 ,0)$ form a copy of the constant 
sheaf (with stalks $\C$) on $X^{an}$. Furthermore there is a continuous section $s_0$ with
\[
s(x) = \left\{ \begin{array}{lll}
0 & \text{if} & x \in X \setminus \{x_0\} \\
(-1,1) & \text{if} & x = x_0 .
\end{array} \right. 
\]
These specific sections provide an isomorphism between $S_B$ and the direct sum of the constant
sheaf on $X^{an}$ and the skyscraper sheaf concentrated at $x_0$. Hence the sheaf cohomology
$\check H^* (X^{an} ; S_B)$ is isomorphic to
\[
\check H^* (X^{an} ; \C) \oplus \check H^* (\{x_0\} ; \C) .
\]
Let us compare this to $HP_* (B)$. We abbreviate $\check H^{ev} = \bigoplus_{m \in \Z} 
\check H^{2m}$ and $\check H^{odd} = \bigoplus_{m \in \Z} \check H^{2m + 1}$. Applying 
Theorem \ref{thm:2.2} and the excision property of $HP_*$ to the extension
\[
0 \to \{ f \in \mc O (X) : f(x_0) = 0 \} \otimes M_2 (\C) \; \to \; B 
\; \to \; \mc O (\{x_0\}) \otimes \C^2 = \C^2 \to 0
\]
we find an exact hexagon
\begin{equation}\label{eq:hexagon}
\hexagon{\check H^{ev}(X,\{x_0\};\C)}{HP_0 (B)}{\check H^{ev}(\{x_0\};\C^2) \cong \C^2}{
\check H^{odd}(\{x_0\};\C^2) = 0}{HP_1 (B)}{\check H^{odd}(X,\{x_0\};\C)}
\end{equation}
By comparing \eqref{eq:hexagon} with the hexagon associated to the algebra extension
\[
0 \to \{ f \in \mc O (X) : f(x_0) = 0 \} \to \mc O (X) \to \mc O (\{x_0\}) = \C \to 0
\]
one sees that the right vertical map in \eqref{eq:hexagon} is zero. We deduce that 
\[
HP_0 (B) \cong \check H^{ev}(X^{an};\C) \oplus \C \quad \text{and} \quad 
HP_1 (B) \cong \check H^{odd}(X^{an};\C) .
\]
The resulting isomorphism between $HP_* (B)$ and the sheaf cohomology of $S_B$ is not quite 
canonical, since it depends on the decomposition of $S_B$ into a constant sheaf and a skyscraper 
sheaf.\\[2mm]

In general, let $\check H^p (X^{an}; S_A )$ denote the $p$-th \u Cech-cohomology group of the 
sheaf $S_A$. Since the space $X^{an}$ is so nice, there is no need to distinguish between the 
\u Cech-cohomology and the sheaf cohomology of $S_A$.

\begin{thm}\label{thm:2.3}
Let $A$ be a finite type $\mc O (X)$-algebra.
There exists a spectral sequence converging to both $HP_n (A)$ and 
$\bigoplus_{m \in \Z} \check H^{n+2m} (X^{an};S_A )$. In particular these 
finite dimensional vector spaces are \textup{(}unnaturally\textup{)} isomorphic.
\end{thm}
\emph{Proof.}
Using the ``abelian filtrations'' of $A$ studied in \cite{KNS}, we can set up and compare two
spectral sequences, one converging to $HP_* (A)$ and the other to $H^* (X^{an} ; S_A )$.
They differ in level zero, but agree on all higher levels. The details can be found in the
author's PhD thesis \cite[Section 2.2]{Sol1}. $\qquad \Box$
\\[2mm]

With Theorem \ref{thm:2.3} one can determine $HP_* (A)$ from Prim$(A)$ by geometric means. 
In particular it becomes unnecessary to consider ideals as in the abelian filtration of $A$, 
which can be hard to determine explicitly.
\\[2mm]
\textbf{Examples.}\\
First the commutative case, where $I$ is an ideal of $\mc O (X)$ and $Y \subset X$
is the zero locus of $I$. Then $S_I$ is the direct image, with respect to the inclusion map 
$X \setminus Y \to X$, of the constant sheaf with fiber $\C$ on $X \setminus Y$.
Thus it appears that we recover Theorem \ref{thm:2.1} from Theorem \ref{thm:2.3}, but in fact the
former is an essential ingredient in the proof of the latter.
\vspace{1mm}

Next suppose that a finite group $G$ acts on the affine variety $X$, and consider the finite type 
$\mc O (X/G)$-algebra $B := \mc O (X) \rtimes G$. Using Theorem \ref{thm:1.1}.e one can construct 
a continuous bijection from the extended quotient $\tilde X / G$ to Prim$(\mc O (X) \rtimes G)$. 
We note that in general this map is neither canonical nor an homeomorphism. 

Consider the constant sheaf with fiber $\C$ on $\tilde X / G$. We claim that the direct image 
$\mc F$ of this sheaf under the natural map $\tilde X / G \to X / G$ is isomorphic to the sheaf $S_B$. 
By Theorem \ref{thm:1.1}.e both have isomorphic stalks over any point of $X/G$.
To construct a sheaf map $\mc F \to S_B$, consider any section $f$ of $\mc F$. For $Gx \in X /G$ 
in the domain of $f$ there is a unique virtual $G_x$-representation
$\pi_x \in R (G_x ) \otimes_\Z \C \cong S_B (Gx)$ with trace $\mr{tr}(\pi_x , g) = f(g,x)$. We define 
a section $s$ of $S_B$ by $s(x) = \pi_x$. Now $f \mapsto s$ is a sheaf map which induces bijections 
on the stalks, so it is an isomorphism $\mc F \to S_B$.
\\[2mm]

For later use we will now calculate the cohomology of certain sheaves of the type $S_A$. 
Let $X$ be a countable, locally finite, finite dimensional CW-complex space endowed with 
an action of a finite group $G$. Let $V$ be a finite
dimensional complex vector space and put $A = C(X) \otimes_\C \mr{End}_\C (V)$.
Suppose furthermore that we have elements $u_g \in A^\times$ such that 
\[
(g \cdot a)(x) := u_g (x) a (g^{-1} x) u_g^{-1} (x) \qquad g \in G ,a \in A, x \in X
\]
defines an action of $G$ on $A$ by algebra automorphisms. In particular, for every $x \in X$ we
get a projective $G_x$-representation $V_x$, by $g \mapsto u_g (x) \in PGL (V)$.

The algebra $A^G$ of $G$-invariant elements in $A$ can be considered as a topological analogue 
of a finite type algebra. Its topological $K$-theory $K_* \big( A^G \big)$ behaves much like the 
\pch \ of finite type algebras. Indeed, we can define a sheaf $S_{A^G}$ over $X/G$ in the same way 
as in Definition \ref{def:SA}. According to \cite[Theorem 2.24]{Sol1} there is an isomorphism
\begin{equation}\label{eq:2.3}
K_n \big( A^G \big) \otimes_\Z \C \cong {\ts \bigoplus_{m \in \Z}} \check H^{n+2m} (X/G ; S_{A^G}) ,
\end{equation}
if at least one side has finite dimension.

\begin{lem}\label{lem:2.4}
Suppose that $X$ is $G$-equivariantly contractible to a point $x_0 \in X$. Then
\[
\check H^n (X/G ; S_{A^G} ) = 0 \qquad \forall n > 0 
\]
and the evaluation at $x_0$ induces a linear bijection
\[
\check H^0 (X/G ; S_{A^G} ) \to S_{A^G} (x_0 ) \cong HH_0 \big( \mr{End}_G (V_{x_0}) \big) .
\]
\end{lem}
\emph{Proof.}
By \cite[Lemma 2.26]{Sol1} the evaluation map $A^G \to \mr{End}_G (V_{x_0})$ is a homotopy 
equivalence. In particular it induces an isomorphism on topological $K$-theory, so by \eqref{eq:2.3}
\[
{\ts \bigoplus_{m \in \Z}} \check H^{n+2m} (X/G ; S_{A^G} ) \cong
K_n \big( A^G \big) \otimes_\Z \C \cong K_n \big( \mr{End}_\C (V_{x_0}) \big) \otimes_\Z \C .
\]
For odd $n$ this is vector space is zero, while for even $n$ it is naturally isomorphic to
\[
HH_0 \big( \mr{End}_G (V_{x_0}) \big) \cong S_{A^G}(x_0 ) .
\]
Since the map
\[
\check H^0 (X/G ; S_{A^G} ) \to \check H^0 (\{ x_0 \} ; S_{A^G} ) = S_{A^G} (x_0 ) 
\]
is already surjective, it must be bijective, and $\check H^n (X/G ; S_{A^G} )$ must be zero
for $n > 0. \qquad \Box$
\vspace{4mm}

\section{Some representation theory}
\label{sec:5}

We discuss the representation theory of $\mh H$ and $\mh H'$, via parabolically induced
representations and intertwining operators. Most results here rely on the author's
previous work \cite{Sol2}. 

As on page \pageref{eq:3.6}, let $\Gamma$ be a finite group of diagram automorphisms of
$(R,\Pi )$. Consider the group $W' = \Gamma \ltimes W$ and the algebra 
$\mh H' = \Gamma \ltimes \mh H$. 

\begin{prop}\label{prop:5.3}
\begin{description}
\item[a)] We have
\[
Z \big( \mh H' \big) \supset S (\mf t^* )^{W'} = \mc O (\mf t / W') ,
\]
with equality if the action of $W'$ on $\mf t^*$ is faithful.
\item[b)] $\mh H'$ is a finite type $\mc O (\mf t / W')$-algebra.
\end{description}
\end{prop}
\emph{Proof.}
a) Lusztig \cite[Theorem 6.5]{Lus1} proved the corresponding statement for the graded 
Hecke algebra $\tilde{\mh H} ( \tilde{\mc R} )$ with formal parameters $\mb k_\alpha$. 
We will adept his argument to our setting. 

On the one hand, it is clear from the multiplication rules \eqref{eq:3.1} and \eqref{eq:3.6} 
that $S (\mf t^* )^{W'} \subset Z (\mh H' )$. On the other hand, consider any
\[
z = {\ts \sum_{w' \in W'}} \, w' p_{w'} \in Z (\mh H') \text{ , with } p_{w'} \in S(\mf t^* ) .
\]
Assume that $W'$ acts faithfully on $\mf t^*$ and that $p_{w'} \neq 0$ for some $w' \neq e$. 
Among these, pick one $w'_0 = \gamma_0 w_0 \in \Gamma \ltimes W$ such that $w_0$ is maximal 
with respect to the Bruhat order $\leq$ for $(W,\Pi )$. For $x \in \mf t^*$ a repeated 
application of \eqref{eq:3.1} yields
\[
\sum_{w' \in W'} w' p_{w'} x = z x = x z = \sum_{w' = \gamma w \in W'} \gamma \Big( w \, w'(x) 
+ \sum_{u \leq w} u f(\gamma (x),w,u) \Big) p_{\gamma w} ,
\]
for suitable polynomials $f(\gamma (x),w,u) \in S (\mf t^* )$. Comparing the coefficients at 
$w'_0$ on both sides, we see that $w'_0 (x) = x$. Since $x \in \mf t^*$ was arbitrary and the
action was faithful, this implies $w'_0 = e$, contrary to our assumption. Therefore 
$z \in Z (\mh H' ) \subset S(\mf t^* )$. But now $w' z = z w'$ for all $w' \in W'$ shows that 
$z \in S (\mf t^* )^{W'}$.\\
b) Clearly $\mh H'$ is an $S(\mf t^* )$-module of finite rank, so it suffices to show that
$S(\mf t^* ) = \mc O (\mf t )$ is of finite rank over $S(\mf t^* )^{W'} = \mc O (\mf t / W' )$.
For any $f \in S(\mf t^* )$ the monic polynomial
\[
P_f (X) := \prod\nolimits_{w \in W'}(X - w(f))
\]
has coefficients in $S(\mf t^* )^{W'}$. Obviously $P_f (f) = 0$, so $S(\mf t^* )$ is integral
over $S(\mf t^* )^{W'}$. Since $S(\mf t^* )$ is finitely generated as an algebra, it follows
that it is an $S(\mf t^* )^{W'}$-module of finite rank. $\qquad \Box$
\\[2mm]

If an $\mh H'$-module $V$ admits an $S (\mf t^* )^{W'}$-character $W' \lambda \in \mf t / W'$,
then we will refer to $W' \lambda$ as the central character of $V$. We say that a central
character is real if it lies in $\mf a / W'$.

A more subtle tool to study $\mh H'$-modules is restriction to the commutative
subalgebra $S (\mf t^* ) \subset \mh H'$. Let $(\pi ,V)$ be an $\mh H'$-module and pick
$\lambda \in \mf t$. The $\lambda$-weight space of $V$ is
\[
V_\lambda = \{ v \in V : \pi (x) v = \inp{x}{\lambda} v \; \forall x \in \mf t^* \} \,.
\] 
We call $\lambda$ an $S (\mf t^* )$-weight of $V$ if $V_\lambda \neq 0$.

Let $P \subset \Pi$ be a set of simple roots. They form a basis of a root subsystem
$R_P \subset R$ with Weyl group $W_P \subset W$. Let $\mf a_P \subset \mf a$ and 
$\mf a_P^* \subset \mf a^*$ be the real spans of respectively $R_P^\vee$ and $R_P$.
We denote the complexifications of these vector spaces by $\mf t_P$ and $\mf t_P^*$, 
and we write 
\[
\begin{array}{lllll}
\mf t^P & = & (\mf t_P^* )^\perp & = &
\{ \lambda \in \mf t : \inp{x}{\lambda} = 0 \: \forall x \in \mf t_P^* \} \,, \\
\mf t^{P*} & = & (\mf t_P )^\perp & = & 
\{ x \in \mf t^* : \inp{x}{\lambda} = 0 \: \forall \lambda \in \mf t_P \} \,.
\end{array}
\]
We define the degenerate root data 
\begin{align}
& \tilde{\mc R}_P = (\mf a_P^* ,R_P ,\mf a_P ,R_P^\vee ,P) \,, \\
& \tilde{\mc R}^P = (\mf a^* ,R_P ,\mf a ,R_P^\vee ,P) \,,
\end{align}
and the graded Hecke algebras
\begin{align}
\mh H_P = \mh H (\tilde{\mc R}_P ,k) \,, \\
\mh H^P = \mh H (\tilde{\mc R}^P ,k) \,.
\end{align}
Notice that the latter decomposes as a tensor product of algebras: 
\begin{equation}
\mh H^P = S (\mf t^{P*} ) \otimes \mh H_P .
\end{equation}
In particular every irreducible $\mh H^P$-module is of the form $\mh C_\lambda \otimes V$, 
where $\lambda \in \mf t^P$ and $V$ is an irreducible $\mh H_P$-module.
In general, for any $\mh H_P$-module $(\rho, V_\rho )$ and $\lambda \in \mf t^P$
we denote the action of $\mh H^P$ on $\mh C_\lambda \otimes V_\rho$ by 
$\rho_\lambda$. We define the parabolically induced module
\begin{equation}
\pi (P,\rho ,\lambda ) = \mr{Ind}_{\mh H^P}^{\mh H} (\mh C_\lambda \otimes V_\rho )
= \mr{Ind}_{\mh H^P}^{\mh H} (\rho_\lambda ) = \mh H \otimes_{\mh H^P} V_{\rho_\lambda} \,.
\end{equation}
Since the complex vector space $\mf t$ has a distinguished real form $\mf a$, we can
decompose any $\lambda \in \mf t$ unambiguously as
\begin{equation}
\lambda = \Re (\lambda ) + i \Im (\lambda ) \text{ with } 
\Re (\lambda ), \Im (\lambda ) \in \mf a \,.
\end{equation}
We define the positive cones 
\begin{equation}
\begin{array}{lll}
\mf a^{*+} & = & \{ x \in \mf a^* : 
  \inp{x}{\alpha^\vee} \geq 0 \: \forall \alpha \in \Pi \} \,, \\
\mf a_P^+ & = &  \{ \mu \in \mf a_P : 
  \inp{\alpha}{\mu} \geq 0 \: \forall \alpha \in P \} \,, \\
\mf a^{P+} & = & \{ \mu \in \mf a^P : 
  \inp{\alpha}{\mu} \geq 0 \: \forall \alpha \in \Pi \setminus P \} \,, \\
\mf a^{P++} & = & \{ \mu \in \mf a^P : 
  \inp{\alpha}{\mu} > 0 \; \forall \alpha \in \Pi \setminus P \} \,.
\end{array}
\end{equation}
The antidual of $\mf a^{*+}$ is
\begin{equation} 
\mf a^- = \{ \lambda \in \mf a : \inp{x}{\lambda} \leq 0 \: \forall x \in \mf a^{*+} \} = 
\big\{ {\ts \sum_{\alpha \in \Pi}} \lambda_\alpha \alpha^\vee : \lambda_\alpha \leq 0 \big\} \,.
\end{equation}
The interior $\mf a^{--}$ of $\mf a^-$ equals
\[
\big\{ {\ts \sum_{\alpha \in \Pi}} \lambda_\alpha \alpha^\vee : \lambda_\alpha < 0 \big\}
\]
if $\Pi$ spans $\mf a^*$, and is empty otherwise.
A finite dimensional $\mh H$-module $V$ is called tempered if $\Re (\lambda ) \in \mf a^-$,
for all weights $\lambda$. More restrictively we say that $V$ belongs to the
discrete series if it is irreducible and $\Re (\lambda ) \in \mf a^{--}$, again for all
weights $\lambda$.

For now on we assume that the root system $R$ is crystallographic and that $k$ is real, that is, 
$k_\alpha \in \R$ for all $\alpha \in \Pi$. 
Of particular importance are the parabolically induced representations 
$\pi (P,\delta ,\lambda )$ where $(\delta ,V_\delta )$ is a discrete series representation 
of $\mh H_P$. We call such a triple $(P,\delta ,\lambda )$ an induction datum and we denote 
the space of these by $\tilde \Xi$. For $\xi ,\eta \in \tilde \Xi$ we write 
$\xi \cong \eta$ if $\xi = (P,\delta,\lambda )$ and $\eta = (P,\sigma ,\lambda)$ with 
$\delta \cong \sigma$ as $\mh H_P$-modules. The space of unitary induction data is
\[
\tilde \Xi_u := \{ (P,\delta ,\lambda ) \in \tilde \Xi : \lambda \in i \mf a^P \} .
\]
There exists a natural involution on $\mh H$ such that $\pi (\xi )$ is a unitary module
if and only if $\xi \in \tilde \Xi_u$, see \cite[Section 6]{Sol2}. Moreover, by 
\cite[Proposition 7.2]{Sol2} every irreducible $\mh H$-module is a quotient of some 
$\pi (P,\delta ,\lambda )$ with $\lambda \in \mf a^{P+}$. Thus discrete series 
representations can be considered as the basic objects for constructing
$\mh H$-modules, and it is certainly useful to know more about them.

\begin{thm}\label{thm:5.1}
There are only finitely many equivalence classes of discrete series representations,
and the $S(\mf t^* )$-weights of a discrete series representation all lie in $\mf a$.
\end{thm}
\emph{Proof.}
See Lemma 2.13 and Corollary 2.14 of \cite{Slo}, or \cite[Theorem 6.2]{Sol2}.
We note that these proofs depend on the classification of possible $S(\mf t^*)$-weights
in \cite[Section 4]{HeOp}. $\qquad \Box$
\\[2mm]

The strategy of parabolic induction can be extended from $\mh H$ to 
$\mh H' = \Gamma \ltimes \mh H$. For an induction datum $(P,\delta ,\lambda) \in \Xi$ we put
\[
\pi' (P,\delta ,\lambda) = \mr{Ind}_{\mh H}^{\mh H'} \pi (P,\delta ,\lambda) =
\mr{Ind}_{\mh H^P}^{\mh H'} (\delta_\lambda ) .
\]
For $P,Q \subset \Pi$ we write
\[
W'(P,Q) = \{ w \in W' : w(P) = Q \} .
\]
Any $w \in W'(P,Q)$ induces algebra isomorphisms
\begin{align}
& \nonumber \psi_w : \mh H_P \to \mh H_Q \\
& \label{eq:5.1} \psi_w : \mh H^P \to \mh H^Q \\
& \nonumber \psi_w (x u) = w(x) \, w u w^{-1} 
\qquad x \in \mf t^* , u \in W_P \,.
\end{align}
We agree to call a $\Gamma \ltimes \mh H$-module tempered if its restriction
to $\mh H$ is tempered. Let $(\pi ,V)$ be any $\mh H$-module and let $\gamma \in \Gamma$. 
Because $\Gamma \cdot \mf a^- = \mf a^-$, the $\mh H$-module $(\pi \circ \psi_\gamma^{-1} ,V)$ 
is tempered, respectively discrete series, if and only if $\pi$ is tempered, respectively
discrete series.

For $w \in W'(P,Q)$ and $(P,\delta ,\lambda ) \in \tilde \Xi$ we have
\[
w (P,\delta ,\lambda ) := (Q ,\delta \circ \psi_w^{-1} ,w(\lambda )) \in \tilde \Xi .
\]
This defines a partial action of $W'$ on $\tilde \Xi$, which preserves $\tilde \Xi_u$.
For $\xi \in \tilde \Xi$ we write
\[
W'_\xi = \{ w \in W' : w (\xi ) \cong \xi \} .
\]
Induction data that are $W'$-associate usually, but not always, yield equivalent $\mh H'$-modules:

\begin{prop}\label{prop:5.2}
Let $w, P, \delta$ and $\lambda$ be as above. There exists an intertwining operator 
\[
\pi' (w, P, \delta ,\lambda ) : \pi' (P,\delta ,\lambda ) \to \pi'(w (P,\delta ,\lambda )) ,
\]
which is rational as a function of $\lambda \in \mf t^P$.
It is regular and invertible for $\lambda$ in a nonempty Zariski-open subset of $\mf t^P$.
\end{prop}
\emph{Proof.}
See Proposition 3.3 and page 28 of \cite{Sol2}. $\qquad \Box$
\\[2mm]
Recall that $R$ was crystallographic and that $k_\alpha \in \R$ for all $\alpha \in \Pi$. 
For unitary induction data much more can be said:

\begin{thm}\label{thm:5.4}
Let $\xi ,\eta \in \tilde \Xi_u$.
\begin{description}
\item[a)] The $\mh H'$-module $\pi' (\xi )$ is tempered and completely reducible.
\item[b)] The operators
\[
\{ \pi' (w,\xi ) : w \in W' , w (\xi ) \cong \eta \}
\]
are well-defined and invertible.
\item[c)] These operators span $\mr{Hom}_{\mh H'} (\pi' (\xi ), \pi' (\eta ))$.
\item[d)] $\{ \pi' (\xi )(h) : h \in \mh H' \} = \{ M \in \mr{End}_\C (\pi' (\xi )) : 
\pi' (w,\xi ) M \pi' (w,\xi )^{-1} = M \; \forall w \in W'_\xi \}$.
\end{description}
\end{thm}
\emph{Proof.}
a) is the combination of Corollary 6.4 and Theorem A.1.c of \cite{Sol2}.\\ 
b) and c) are in \cite[Theorem 8.2.c]{Sol2}.\\
d) By parts b) and c) the right-hand side is the bicommutant of the left-hand side.
By part a) $\pi' (\xi) (\mh H' )$ is a semisimple algebra, hence equal to its own
bicommutant in $\mr{End}_\C (\pi' (\xi )). \qquad \Box$
\vspace{4mm}

\section{The primitive ideal spectrum}
\label{sec:6}

We will describe the primitive ideal spectrum of $\mh H' = \Gamma \rtimes \mh H$ as explicitly
as possible. First we stratify it and we describe the strata as topological spaces. After that
we will determine the cohomology of the sheaf corresponding to Prim$(\mh H')$ and draw some
conclusions about the \pch \ of $\mh H'$. Throughout this section we assume that the root system 
$R$ is crystallographic and that $k_\alpha \in \R$ for all $\alpha \in \Pi$.

We denote the central character of an (irreducible) $\mh H_P$-representation by
$cc_P (\delta ) \in \mf t_P / W_P$. Since $W_P$ acts orthogonally on $\mf t$, all 
$S (\mf t^*_P )$-weights of $\delta$ have the same norm, which allows us to write
$\| cc_P (\delta ) \|$ without ambiguity.

\begin{thm}\label{thm:6.1} 
\textup{\cite[Theorem 8.3.b]{Sol2}} \\
Let $\rho$ be an irreducible $\mh H'$-module. There exists a unique $W'$-association class
$W' \xi_\rho = W' (P_\rho,\delta_\rho ,\lambda_\rho )$ in $\tilde \Xi / W'$ such that:
\begin{itemize}
\item $\rho$ is a constituent of $\pi' (\xi_\rho )$,
\item $\| cc_{P_\rho} (\delta_\rho ) \|$ is maximal among induction data with the first property.
\end{itemize}
\end{thm}

\begin{lem}\label{lem:6.3}
Let $\rho$ be an irreducible $\mh H'$-module with $W' \xi_\rho = W' (P_\rho ,\delta_\rho ,\lambda_\rho )$.
\begin{description}
\item[a)] $\rho$ is tempered if and only if $\xi_\rho \in \tilde \Xi_u$, which is equivalent to
$\lambda_\rho \in i \mf a$.
\item[b)] The central character of $\rho$ lies in $\mf a / W'$ if and only if $\lambda_\rho \in \mf a$.
\end{description}
\end{lem}
\emph{Proof.}
a) follows from Proposition 7.3.c and Theorem 8.3.a of \cite{Sol2}.\\
b) The central character of $\pi' (P,\delta ,\lambda )$ is $W' (cc_P (\delta ) + \lambda )$
and by Theorem \ref{thm:5.1} $cc_P (\delta ) \in \mf a / W_P . \qquad \Box$
\\[2mm]

For every $P \subset \Pi$, every discrete series representation $\delta$ of $\mh H_P$ and every
$U \subset \mf t^P$ we get a subset of Prim$(\mh H')$:
\[
\mr{Prim}_{(P,\delta ,U)} (\mh H' ) := \{ \ker \rho \in \mr{Prim}(\mh H' ) : 
W' \xi_\rho \cap (P,\delta ,U ) \neq \emptyset \} .
\]
For $U = \mf t^P$ or $U = \{ \lambda \}$ we abbreviate this to 
$\mr{Prim}_{P,\delta} (\mh H' )$ or $\mr{Prim}_{(P,\delta ,\lambda)} (\mh H' )$.

The $\mr{Prim}_{P,\delta} (\mh H' )$ form a partition of Prim$(\mh H' )$, but they are not closed. 
However, the boundary of $\mr{Prim}_{P,\delta} (\mh H' )$ can only contain elements of 
$\mr{Prim}_{Q,\sigma} (\mh H' )$ if $\| cc_Q (\sigma ) \| > \| cc_P (\delta ) \|$. 
By \cite[Theorem 6.2.b]{Sol2} there are only finitely many $W'$-association classes of pairs 
$(P,\delta )$. Pick one $(P_i ,\delta_i )$ from every $W'$-association class, and order them such that 
\[
\| cc_{P_i} (\delta_i ) \| \leq \| cc_{P_j} (\delta_j ) \| \text{ if } j \leq i . 
\]
Then $\bigcup_{j \leq i} \mr{Prim}_{P_j ,\delta_j} (\mh H' )$ is closed in Prim$(\mh H')$
for every $i$. This defines a stratification of Prim$(\mh H' )$, with strata 
$\mr{Prim}_{P,\delta} (\mh H' )$. This stratification is the analogue of a stratification of the smooth
dual of a reductive $p$-adic group \cite[Lemma 2.17]{Sol3}, which for $GL_n$ can already be
found in \cite[Section 9]{ScZi}.

The group
\[
W'_\delta := \{ w \in W' : w(P) = P , \delta \circ \psi_w^{-1} \cong \delta \}
\]
acts on $\mf t^P$ and $\lambda$'s in the same $W'_\delta$-orbit give rise to the same elements of
$\mr{Prim}_{P,\delta} (\mh H' )$. Hence there is a canonical map
\[
\mr{Prim}_{P,\delta}(\mh H' ) \to \mf t^P /W'_\delta ,
\]
which is continuous, surjective and finite-to-one.

\begin{prop}\label{prop:6.2} 
\textup{\cite[Proposition 10.1]{Sol2}} \\
Suppose that $U \subset \mf t^P$ satisfies
\begin{itemize}
\item $U \to \{-1,0,1\} : \lambda \mapsto \mr{sign} \inp{\Re (\lambda )}{\alpha}$
is constant for all $\alpha \in \Pi \setminus P$,
\item every $\lambda \in U$ has the same stabilizer in $W'_\delta$.
\end{itemize}
Then $\mr{Prim}_{(P,\delta ,U)} (\mh H' )$ is homeomorphic to $U / W'_\delta \times 
\mr{Prim}_{(P,\delta ,\lambda_0 )} (\mh H' )$, for any $\lambda_0 \in U$.
\end{prop}

Thus $\mr{Prim}_{P,\delta}(\mh H' )$ looks like the quotient of the vector space $\mf t^P$ by
$W'_\delta$, but with certain linear subspaces carrying a multiplicity.

Consider the sheaf $S_{\mh H'}$ over $\mf t / W'$, as in Definition \ref{def:SA}.
It has a subsheaf $S_{P,\delta}$ consisting of all those sections of $S_{\mh H'}$ that use only 
elements of $\mr{Prim}_{P,\delta}(\mh H' )$. Alternatively we can obtain $S_{P,\delta}$ as the
sheaf corresponding to a certain subquotient algebra of $\mh H'$, but there is no need to use this 
subquotient.

\begin{lem}\label{lem:6.4}
\[
\check H^n (\mf t / W' ; S_{P,\delta} ) = 0 \qquad \forall n > 0 ,
\]
and evaluation at $\mu = W' cc_P (\delta ) \in \mf t /W'$ induces a linear bijection
\[
\check H^0 (\mf t / W' ; S_{P,\delta} ) \longrightarrow S_{P,\delta} (\mu ) =
\C \{ I \in \mr{Prim}(\mh H' ) : I \supset \ker (\pi' (P,\delta ,0)) \} .
\]
\end{lem}
\emph{Proof.}
By construction $S_{P,\delta}$ is the direct image of a sheaf $\mc F_\delta$ on 
$\mf t^P / W'_\delta$, with respect to the map
\begin{equation}\label{eq:6.1}
\mf t^P / W'_\delta \to \mf t^P : W'_\delta \lambda \mapsto W' (\lambda + cc_P (\delta )) .
\end{equation}
According to \cite[Theorem 11.1]{Bre} this induces a natural isomorphism
\begin{equation}\label{eq:6.2}
\check H^* (\mf t /W' ; S_{P,\delta} ) \cong \check H^* \big( \mf t^P /W'_\delta ; \mc F_\delta \big) .
\end{equation}
By Proposition \ref{prop:5.2} and Theorem \ref{thm:5.4}.b we can find a ball $B_r \subset \mf t^P$
around 0, of radius $r > 0$, such that all the operators $\pi'(w,P,\delta ,\lambda )$ with 
$w \in W'_\delta$ and $\lambda \in B_r$ are regular and invertible. By Proposition \ref{prop:6.2}
the homeomorphism
\[
B_r \to \mf t^P : \lambda \mapsto \tan \Big( \frac{\pi \| \lambda \|}{2 r} \Big) \lambda
\]
induces a sheaf isomorphism $\mc F_\delta |_{B_r /W'_\delta} \to \mc F_\delta$. Hence
\begin{equation}
\check H^* (\mf t /W' ; S_{P,\delta} ) \cong \check H^* (B_r /W'_\delta ; \mc F_\delta ) .
\end{equation}
Let $V = \mh H' \otimes_{\mh H^P} V_\delta$ be the vector space on which all the 
representations $\pi'(P,\delta ,\lambda )$ are realized, and consider the Fr\'echet algebra 
$A = C(B_r ) \otimes \mr{End}_\C (V)$. For $w \in W'_\delta$ we define 
\[
u_w \in C(B_r ) \otimes \mr{Aut}_\C (V) \quad \mr{by} \quad
u_w (\lambda ) = \pi' (w,P,\delta ,w^{-1} \lambda ) .
\]
Now $W'_\delta$ acts on $A$ by 
\[
(w \cdot a)(\lambda ) = \pi' (w,P,\delta ,w^{-1} \lambda ) a(w^{-1} \lambda ) 
\pi' (w,P,\delta ,w^{-1} \lambda )^{-1} ,
\]
which can be abbreviated to $w \cdot a = u_w (a \circ w^{-1} ) u_w^{-1}$. Thus we are in the 
setting of page \pageref{eq:2.3}. We will relate our sheaves to the algebra $A' = A^{W'_\delta}$ 
of $W'_\delta$-invariants. By Theorem \ref{thm:5.4}.d we have 
\[
\{ \pi' (P,\delta ,\lambda )(h) : h \in \mh H' \} = \big\{ a (\lambda ) : a \in A' \big\} 
\]
for all $\lambda \in i \mf a^P$. Hence the sheaves $F_\delta$ and $S_{A'}$ have the 
same restriction to $(B_r \cap i \mf a^P ) / W'_\delta$. For every $\lambda \in B_r$ the dimension
of $S_{A'} (W'_\delta \lambda )$ equals the number of inequivalent constituents of the
projective $(W'_\delta )_\lambda$-representation $w \mapsto u_w (\lambda )$. Since $(\mf t^P )^G$
is connected for every $G \subset W'_\delta$, this number depends only on $\lambda$ via the 
isotropy group $(W'_\delta )_\lambda$. Combining this with Proposition \ref{prop:6.2}, we see
that the sheaves $\mc F_\delta |_{B_r / W'_\delta}$ and $S_{A'}$ are isomorphic. 
According to Lemma \ref{lem:2.4} 
\[
\check H^n \big( B_r / W'_\delta ; S_{A'} \big) = 0 \qquad \forall n > 0
\]
and the evaluation map 
\[
\check H^0 \big( B_r / W'_\delta ; S_{A'} \big) \to S_{A'} (0)
\]
is bijection. Consequently
\[
\check H^n (\mf t /W' ; S_{P,\delta} ) =  \check H^n (B_r /W'_\delta ; \mc F_\delta ) 
\qquad \forall n > 0
\]
and the maps
\[
\check H^0 (\mf t^P /W' ; \mc F_\delta ) =  \check H^0 (B_r /W'_\delta ; \mc F_\delta ) 
\to \mc F_\delta (0)
\]
are isomorphisms. In view of \eqref{eq:6.1} and \eqref{eq:6.2} this shows that the 
corresponding evaluation
\[
\check H^0 (\mf t /W' ; S_{P,\delta} ) \to S_{P,\delta} (W' cc_P (\delta ))
\]
is also bijective. $\qquad \Box$
\\[2mm]

We denote the finite set of equivalence classes of irreducible tempered $\mh H'$-modules
with real central character by $\mr{Irr}_0 (\mh H' )$. In every such equivalence class we pick
one $\mh H'$-module $(\pi ,V)$ and we confuse $\mr{Irr}_0 (\mh H' )$ with this set of modules.
Let $\C \, \mr{Irr}_0 (\mh H' )$ be the complex vector space with basis $\mr{Irr}_0 (\mh H' )$.

By Theorem \ref{thm:5.4}.a and Lemma \ref{lem:6.3} $\mr{Irr}_0 (\mh H' )$ consists precisely of
the irreducible summands of the modules $\pi' (P,\delta ,0)$, where $P \subset \Pi$ and $\delta$ 
is a discrete series representation of $\mh H_P$. Notice that all the primitive ideals appearing 
in Lemma \ref{lem:6.4} correspond to such modules.

\begin{thm}\label{thm:6.5} 
As before we assume that $R$ is crystallographic and that $k$ is real.
\begin{description}
\item[a)] $\check H^n ( \mf t / W' ; S_{\mh H'} ) = 0$ for all $n > 0$, 
and there exists a canonical linear bijection
\[
\check H^0 (\mf t / W' ; S_{\mh H'} ) \to \C \, \mr{Irr}_0 (\mh H' ) .
\]
\item[b)] The natural algebra homomorphisms 
\[
\C [W'] \longrightarrow \mh H' \longrightarrow 
{\ts \bigoplus_{(\pi ,V) \in \mr{Irr}_0 (\mh H' ) }} \mr{End}_\C (V )
\]
induce isomorphisms on \pch .
\item[c)] The $W'$-representations $\{ V |_{W'} : (\pi ,V) \in \mr{Irr}_0 (\mh H' ) \}$ 
form a $\Q$-basis of the representation ring $R(W' ) \otimes_\Z \Q$.
\end{description}
\end{thm}
\emph{Proof.}
a) On page \pageref{prop:6.2} we defined a stratification of Prim$(\mh H' )$ with terms
\[
F_i := {\ts \bigcup_{j \leq i} } \mr{Prim}_{P_j ,\delta_j} (\mh H' ) .
\]
These give rise to subsheaves $S_i$ of $S_{\mh H'}$, consisting of all sections of $S_{\mh H'}$
that use only elements of $F_i$. By construction $S_i / S_{i-1}$ is isomorphic to the sheaf
$S_{P_i ,\delta_i }$ studied in Lemma \ref{lem:6.4}. This filtration of $S_{\mh H'}$ yields a
spectral sequence that converges to $\check H^* ( \mf t / W' ; S_{\mh H'} )$ and has
initial terms $E^1_{i,j} = \check H^{i+j} \big( \mf t / W' ; S_{P_i ,\delta_i } \big)$.
By Lemma \ref{lem:6.4} $E^1_{i,j} = 0$ unless $i+j = 0$. Hence all the boundary maps 
$\delta^r_{i,j} : E^r_{i,j} \to E^r_{i-r,j+r-1}$ are 0 for $r \geq 1$, and the spectral
sequence degenerates at $E^1_{*,*}$. Lemma \ref{lem:6.4} tells us that
\[
\begin{split}
\check H^0 (\mf t / W' ; S_{\mh H'} ) & \cong {\ts \bigoplus_i} \check H^0 
(\mf t / W' ; S_{P_i ,\delta_i} ) \\
& \cong {\ts \bigoplus_i } \, \C \{ I \in \mr{Prim}(\mh H' ) : 
I \supset \ker (\pi' (P_i,\delta_i ,0)) \} \cong \C \, \mr{Irr}_0 (\mh H' ) ,
\end{split}
\]
where the composite isomorphism is given by evaluating a global section of $S_{\mh H'}$ at all
points corresponding to $\mr{Irr}_0 (\mh H' )$.\\
b) follows from a) and Theorems \ref{thm:2.3} and \ref{thm:4.3}.\\
c) Since $\C [W']$ is finite dimensional and semisimple, 
\[
HP_* \big( \C [W'] \big) \cong HH_0 \big( \C [W'] \big) = \C [W'] / \big[ \C [W'],\C [W'] \big] 
\cong Z \big( \C [W'] \big) \cong R(W' ) \otimes_\Z \C .
\]
Under these isomorphisms an irreducible $W'$-representation $\rho$ corresponds to the class
of the central idempotent $e_\rho \in Z\big( \C [W' ] \big)$ in $HH_0 \big( \C [W'] \big)$. 
Similarly
\[
\begin{split}
HP_* \big( {\ts \bigoplus_{(\pi ,V) \in \mr{Irr}_0 (\mh H' ) } } \mr{End}_\C (V ) \big) & =
{\ts \bigoplus_{(\pi ,V) \in \mr{Irr}_0 (\mh H' ) } } HH_0 ( \mr{End}_\C (V ) ) \\
& \cong {\ts \bigoplus_{(\pi ,V) \in \mr{Irr}_0 (\mh H' ) } } Z ( \mr{End}_\C (V ) ) \; \cong \;
\C \, \mr{Irr}_0 (\mh H' ) .
\end{split}
\]
Let $HP_* (\phi )$ be the map induced by the algebra homomorphism
\[
\phi : \C [W'] \to {\ts \bigoplus_{(\pi ,V) \in \mr{Irr}_0 (\mh H' ) } } \mr{End}_\C (V ) .
\]
By construction
\[
HP_* (\phi ) (e_\rho ) = {\ts \sum_{(\pi ,V) \in \mr{Irr}_0 (\mh H' ) } } 
  \mr{tr}(\pi (e_\rho )) \; (\pi ,V)
\]
By part b) this is a bijection
\[
R(W') \otimes_\Z \C \cong Z \big( \C [W'] \big) \to \C \, \mr{Irr}_0 (\mh H' ) .
\]
We declare the canonical bases $\widehat{W'}$ and $\mr{Irr}_0 (\mh H' )$ of these vector spaces
to be orthonormal. Then the adjoint map
\[
HP_* (\phi )^* : \C \, \mr{Irr}_0 (\mh H' ) \to R(W') \otimes_\Z \C
\]
is also bijective, and it sends $(\pi ,V)$ to $\pi \circ \phi \in R(W' )$. Hence 
\[
\big\{ \pi \circ \phi : (\pi ,V) \in \mr{Irr}_0 (\mh H' ) \big\}
\]
is a $\C$-basis of $R(W') \otimes_\Z \C$ and a $\Q$-basis of $R(W') \otimes_\Z \Q . 
\qquad \Box$


\vspace{4mm}




\begin{thebibliography}{99}
\addcontentsline{toc}{section}{References}


\bibitem[ABP1]{ABP} A.-M. Aubert, P.F. Baum, R.J. Plymen,  
``The Hecke algebra of a reductive $p$-adic group:
a view from noncommutative geometry'',
pp. 1--34 in: \emph{Noncommutative geometry and number theory},
Aspects of Mathematics \textbf{E37},
Vieweg Verlag, 2006

\bibitem[ABP2]{ABP2} A.-M. Aubert, P.F. Baum, R.J. Plymen,
``Geometric structure in the representation theory of $p$-adic groups'',
C. R. Acad. Sci. Paris \textbf{345} (2007), 573--578

\bibitem[BaNi]{BaNi} P.F. Baum, V. Nistor,        
``Periodic cyclic homology of Iwahori--Hecke algebras'',
K-Theory \textbf{27.4} (2002), 329--357

\bibitem[Bre]{Bre} G.E. Bredon.
\emph{Sheaf Theory 2nd ed.},
Graduate Texts in Mathematics \textbf{170},
Springer, 1997

\bibitem[Bry]{Bry} J.L. Brylinski,          
``Cyclic homology and equivariant theories'',
Ann. Inst. Fourier \textbf{37.4} (1987), 15--28

\bibitem[CaEi]{CaEi} H. Cartan, S. Eilenberg,         
\emph{Homological algebra},
Princeton Landmarks in Mathematics,
Princeton University Press, 1956

\bibitem[Cli]{Cli} A.H. Clifford,           
``Representations induced in an invariant subgroup'',
Ann. of Math. \textbf{38} (1937), 533--550

\bibitem[Ciu]{Ciu} D. Ciubotaru,
``On unitary unipotent representations of $p$-adic groups
and affine Hecke algebras with unequal parameters'',
Representation Theory 12 (2008), 453--498

\bibitem[CuRe]{CuRe} C.W. Curtis, I. Reiner,
\emph{Methods of representation theory Vol. I. 
With applications to finite groups and orders},
Pure and Applied Mathematics,
John Wiley \& Sons, New York, 1981

\bibitem[GeJo]{GeJo} E. Getzler, J.D.S. Jones,
``The cyclic homology of crossed product algebras'',
J. reine angew. Math. \textbf{445} (1993), 161--174 

\bibitem[HeOp]{HeOp} G. Heckman, E. Opdam,
``Yang's system of particles and Hecke algebras'',
Ann. of Math. \textbf{145.1} (1997), 139--173

\bibitem[How]{How} R.B. Howlett,
``Normalizers of parabolic subgroups of reflection groups'',
J. London Math. Soc. (2) \textbf{21} (1980), 62--80

\bibitem[Hum]{Hum} J.E. Humphreys,
\emph{Reflection groups and Coxeter groups},
Cambridge Studies in Advanced Mathematics \textbf{29},
Cambridge University Press, 1990

\bibitem[KNS]{KNS} D. Kazhdan, V. Nistor, P. Schneider,        
``Hochschild and cyclic homology of finite type algebras'',
Sel. Math. New Ser. \textbf{4.2} (1998), 321--359

\bibitem[Lod]{Lod} J.-L. Loday,                 
\emph{Cyclic homology 2nd ed.},
Mathematischen Wissenschaften \textbf{301}
Springer Verlag, 1997

\bibitem[Lus1]{Lus0} G. Lusztig,
``Intersection cohomology complexes on a reductive group'',
Invent. Math. \textbf{75.2} (1984), 205--272

\bibitem[Lus2]{Lus1} G. Lusztig,
``Cuspidal local systems and graded Hecke algebras I'',
Publ. Math. Inst. Hautes \'Etudes Sci. \textbf{67} (1988), 145--202

\bibitem[Lus3]{Lus2} G. Lusztig,
``Affine Hecke algebras and their graded version'',
J. Amer. Math. Soc \textbf{2.3} (1989), 599--635

\bibitem[Nis]{Nis} V. Nistor,
``A non-commutative geometry approach to the representation 
theory of reductive $p$-adic groups: Homology of Hecke algebras, 
a survey and some new results'', pp. 301--323 in:
\emph{Noncommutative geometry and number theory},
Aspects of Mathematics \textbf{E37},
Vieweg Verlag, 2006

\bibitem[ScZi]{ScZi} P. Schneider, E.-W. Zink,
``K-types for the tempered components of a $p$-adic general 
linear group. With an appendix by Schneider and U. Stuhler'',
J. Reine Angew. Math. \textbf{517} (1999), 161--208

\bibitem[Slo]{Slo} K. Slooten,
``Generalized Springer correspondence and Green functions 
for type B/C graded Hecke algebras'',
Adv. Math. \textbf{203} (2005), 34--108

\bibitem[Sol1]{Sol1} M.S. Solleveld,
\emph{Periodic cyclic homology of affine Hecke algebras},
Ph.D. Thesis, Universiteit van Amsterdam, 2007, arXiv:0910.1606

\bibitem[Sol2]{Sol3} M.S. Solleveld,
``Periodic cyclic homology of reductive p-adic groups'',
J. Noncommutative Geometry \textbf{3.4} (2009), 501--558

\bibitem[Sol3]{Sol2} M.S. Solleveld
``Parabolically induced representations of graded Hecke algebras'',
arXiv:0804.0433, 2008



\end{thebibliography}
\end{document}